\documentclass[12pt]{amsart}
\usepackage{amsmath}
\usepackage{amsfonts}
\usepackage{amssymb}
\usepackage[all]{xy}
\usepackage{bbding}
\usepackage{txfonts}
\usepackage{amscd}
\usepackage{xspace}
\usepackage[shortlabels]{enumitem}
\usepackage{ifpdf}
\ifpdf
  \usepackage[colorlinks,final,backref=page,hyperindex]{hyperref}
\else
  \usepackage[colorlinks,final,backref=page,hyperindex,hypertex]{hyperref}
\fi
\usepackage{tikz}
\usepackage[active]{srcltx}

\newtheorem{thm}{Theorem}[section]
\newtheorem{lem}[thm]{Lemma}
\newtheorem{cor}[thm]{Corollary}
\newtheorem{pro}[thm]{Proposition}
\newtheorem{ex}[thm]{Example}
\newtheorem{rmk}[thm]{Remark}
\newtheorem{defi}[thm]{Definition}

\newcommand {\emptycomment}[1]{}
\newcommand{\lon }{\,\rightarrow\,}
\newcommand{\be }{\begin{equation}}
\newcommand{\ee }{\end{equation}}
\newcommand{\K}{\mathbb{K}}

\newcommand{\pf}{\noindent{\bf Proof.}\ }
\newcommand {\yh}[1]{{\marginpar{*}\scriptsize\textcolor{purple}{yh: #1}}}

\newcommand {\jf}[1]{{\marginpar{*}\scriptsize\textcolor{blue}{jf: #1}}}

\newcommand{\huaB}{\mathcal{B}}

\newcommand{\huaA}{\mathcal{A}}
\newcommand{\huaL}{\mathcal{L}}
\newcommand{\huaR}{\mathcal{R}}

\newcommand{\huaD}{\mathcal{D}}

\newcommand{\huaO}{\mathcal{O}}

\newcommand{\frad}{\mathfrak {ad}}

\newcommand{\frka}{\mathfrak a}

\newcommand{\frkd}{\mathfrak d}

\newcommand{\g}{\mathfrak g}

\newcommand{\frks}{\mathfrak s}

\newcommand{\frkB}{\mathfrak B}

\newcommand{\frkL}{\mathfrak L}

\newcommand{\frkR}{\mathfrak R}

\def\qed{\hfill ~\vrule height6pt width6pt depth0pt}

\newcommand{\half}{\frac{1}{2}}

\newcommand{\Courant}[1]{\left\llbracket  #1\right\rrbracket }
\newcommand{\Poisson}[1]{\{ #1\}}

\newcommand{\br}[1]{   [ \cdot,    \cdot  ]   }
\newcommand{\id}{\mathrm{id}}

\newcommand{\Hom}{\mathrm{Hom}}

\newcommand{\CYBE}{\mathrm{CYBE}}
\newcommand{\AYBE}{\mathrm{AYBE}}
\newcommand{\PYBE}{\mathrm{PYBE}}

\newcommand{\gl}{\mathfrak {gl}}

\newcommand{\ad}{\mathrm{ad}}

\begin{document}

\title{
{Noncommutative Poisson bialgebras}
 }\vspace{2mm}
\author{Jiefeng Liu, Chengming Bai and Yunhe Sheng
}

\date{}

\maketitle

\begin{abstract}
 In this paper, we introduce the notion of a noncommutative Poisson bialgebra, and establish the equivalence between matched pairs, Manin triples and noncommutative Poisson bialgebras. Using quasi-representations and the corresponding cohomology theory of noncommutative Poisson algebras, we study coboundary noncommutative Poisson bialgebras which leads to
the introduction of the Poisson Yang-Baxter equation. A
skew-symmetric solution of the Poisson Yang-Baxter equation
naturally gives a (coboundary) noncommutative Poisson bialgebra.
Rota-Baxter operators, more generally $\huaO$-operators on
noncommutative Poisson algebras, and noncommutative pre-Poisson
algebras are introduced, by which we construct skew-symmetric
solutions of the Poisson Yang-Baxter equation
in some special noncommutative Poisson algebras obtained from
these structures.
\end{abstract}

\tableofcontents
\section{Introduction}

 This paper aims to study the bialgebra theory for noncommutative Poisson algebras, in particular coboundary ones. Skew-symmetric solutions of the Poisson Yang-Baxter  equation  in certain noncommutative Poisson algebras are constructed using $\huaO$-operators and noncommutative pre-Poisson algebras, which give coboundary noncommutative Poisson bialgebras.

\subsection{Noncommutative Poisson algebras and pre-Poisson algebras}

The notion of a noncommutative Poisson algebra was first given by Xu in \cite{Xu}, which is especially suitable for geometric situations.
\begin{defi}
A {\bf noncommutative Poisson algebra} is a triple $(P,\cdot_P,\{-,-\}_P)$, where $(P,\cdot_P)$ is an associative algebra (not necessarily commutative) and $(P,\{-,-\}_P)$ is a Lie algebra, such that  the Leibniz rule holds:
$$\{x,y\cdot_P z\}_P=\{x,y\}_P\cdot_P z+y\cdot_P\{x,z\}_P,\quad \forall~ x,y,z\in P.$$
\end{defi}
 In \cite{FGV}, Flato, Gerstenhaber and Voronov introduced a more general notion of a Leibniz pair and study its cohomology and deformation theory. In particular, they gave the cohomology theory of a noncommutative Poisson algebra associated to a representation using an innovative bicomplex.  Recently, Bao and Ye developed the cohomology theory of noncommutative Poisson algebras associated to  quasi-representations through Yoneda-Ext groups and projective resolutions in \cite{CohomologyPA1,CohomologyPA2}. Noncommutative Poisson algebras had been studied by many authors from different aspects \cite{Kubo1,Kubo3,Kubo4,YYY}. A Poisson algebra in the usual sense is the one where the associative multiplication on $P$ is commutative. Note that there is another noncommutative analogue of Poisson algebras, namely double Poisson algebras (\cite{Van}), which will not be considered in this paper.

Aguiar introduced the notion of a pre-Poisson algebra in \cite{A2} and constructed many examples. A pre-Poisson algebra contains a Zinbiel algebra and a pre-Lie algebra such that some compatibility conditions are satisfied. Zinbiel algebras, which are also called dual Leibniz algebras, were introduced by Loday in \cite{Lod1}, and further studied in  \cite{Liv,Lod2}. Pre-Lie
algebras are a class of nonassociative algebras coming from the
study of convex homogeneous cones, affine manifolds and affine
structures on Lie groups, and  cohomologies of
associative algebras.  They also appeared in many fields in
mathematics and mathematical physics, such as complex and
symplectic structures on Lie groups and Lie algebras, integrable
systems, Poisson brackets and infinite dimensional Lie algebras,
vertex algebras, quantum field theory and operads. See the survey
\cite{Pre-lie} and the references therein for more details.
 A pre-Poisson algebra gives rise to a  Poisson algebra naturally through the  sub-adjacent commutative associative algebra of the Zinbiel algebra and the sub-adjacent  Lie algebra of the pre-Lie algebra. Conversely,  a Rota-Baxter operator action (more generally an $\huaO$-operator action) on a  Poisson algebra gives rise to a pre-Poisson algebra. We can summarize these relations by the following diagram:
$$
\xymatrix{
\ar[rr] \mbox{{\bf Zinbiel algebra} + pre-Lie algebra }\ar[d]_{\mbox{sub-adjacent~~}}
                && \mbox{Pre-Poisson algebra}\ar[d]_{\mbox{sub-adjacent~~}}\\
\ar[rr] \mbox{{\bf comm associative algebra} + Lie algebra }\ar@<-1ex>[u]_{\mbox{~~Rota-Baxter ~action}}
                && \mbox{ Poisson algebra. }\ar@<-1ex>[u]_{\mbox{~~Rota-Baxter~ action}}}
$$

 \subsection{The bialgebra theory for noncommutative Poisson algebras}

 For a given algebraic structure determined by a set of
multiplications, a bialgebra structure
on this algebra is obtained by a corresponding set of
comultiplications together with a set of compatibility conditions
between the multiplications and comultiplications.  For a finite
dimensional vector space $V$ with the given algebraic structure,
this can be achieved by equipping the dual space
$V^*$ with the same algebraic structure and a set of compatibility
conditions between the structures on $V$ and those on $V^*$.

The great importance of the bialgebra theory and the noncommutative Poisson algebra serves as the main motivation for
our interest in a suitable bialgebra theory for the noncommutative Poisson algebra in this paper.

A good compatibility condition  in a bialgebra is prescribed on the one hand by a strong motivation and potential applications, and on the other
hand by a rich structure theory and effective constructions. In the associative algebra context, an antisymmetric infinitesimal bialgebra~\cite{A1,A3,A4,Bai1}   has the same associative multiplications on $A$ and $A^*$, and the comultiplication being a 1-cocycle on $A$ with   coefficients in the tensor representation $A\otimes A$.  In the Lie algebra context, a   Lie bialgebra  consists of a Lie algebra
$(\g,[-,-]_\g$) and a Lie coalgebra $(\g,\delta)$, where $\delta:\g\to\otimes^2\g$ is a Lie comultiplication, such that the Lie comultiplication being a 1-cocycle on $\g$ with  coefficients in the tensor representation $\g\otimes \g$. See \cite{Drinfeld,EtiSh} for more details about Lie bialgebras and applications in mathematical physics. Thus, the representation theory and the cohomology theory  usually  play  essential roles in the study of a bialgebra theory.

In fact, there has been a bialgebra theory for the usual
(commutative) Poisson algebras, the so-called Poisson bialgebras
(\cite{Poissonbialgebra}), in terms of the representation theory
of Passion algebras.   However, a direct generalization is
not available for the noncommutative Poisson algebras.  In this
paper, we apply  quasi-representations   instead of
representations and the corresponding cohomology theory
(\cite{CohomologyPA1,YYY}) to study noncommutative Poisson
bialgebras. Even though both Lie algebras and associative algebras
admit tensor representations as mentioned above, the tensor
product of two representations of a noncommutative Poisson algebra
is not a representation anymore, but a quasi-representation.  This
is the reason why quasi-representations, not representations, are
the main ingredient  in our study of noncommutative Poisson
bialgebras. On the other hand,  the dual  of the regular
representation $(L,R,\ad)$ of a noncommutative Poisson algebra is
also usually not a representation, but a quasi-representation. We
introduce the concept of a coherent noncommutative Poisson
algebra to overcome this problem. Note that the usual Poisson
algebras are coherent. Thus, the bialgebra theory for
noncommutative Poisson algebras established in this paper contains
all the results of Poisson bialgebras given in
\cite{Poissonbialgebra}.

Moreover,  like the case of Poisson bialgebras, the study of coboundary noncommutative Poisson bialgebras leads to
the introduction of the Poisson Yang-Baxter equation in a coherent noncommutative Poisson algebra.
A skew-symmetric solution of the  {Poisson} Yang-Baxter equation naturally gives a (coboundary) noncommutative Poisson bialgebra.

In addition, coherent noncommutative Poisson algebras are closely related to compatible Lie algebras which play important roles in several fields in mathematics and mathematical physics
(\cite{GS1,GS2,GS3,OS}). In fact, a coherent noncommutative Poisson algebra naturally gives a compatible Lie algebra. Consequently, a noncommutative Poisson bialgebra  gives rise to a compatible Lie bialgebra (\cite{WB}), and there is also a similar relationship in the coboundary cases.

 \subsection{Noncommutative pre-Poisson algebras, Rota-Baxter operators and $\huaO$-operators}
The notion of a dendriform algebra  was introduced by Loday in \cite{Lod1} with motivation from periodicity of algebraic K-theory and operads.

\begin{defi} A {\bf dendriform algebra} is a vector space $A$ with two bilinear maps $\succ:A\otimes A\longrightarrow A$ and $\prec:A\otimes A\longrightarrow A$ such that for all $x,y,z\in A$, the following equalities hold:
  \begin{eqnarray*}
    (x\prec y)\prec z=x\prec (y\succ z+y\prec z),~~
    (x\succ y)\prec z=x\succ(y\prec z),~~
    x\succ(y\succ z)=(x\succ y+x\prec y)\succ z.
  \end{eqnarray*}
  \end{defi}
One can obtain an associative algebra  as well as  a pre-Lie algebra from a dendriform algebra. The relations among dendriform algebras, associative algebras, pre-Lie algebras and Lie algebras are given as follows:
$$
\xymatrix{
 \ar[rr]^{\quad x\ast y=x\succ y-y\prec x}\mbox{dendriform algebra}~(A,\succ,\prec) \ar[d]^{x\cdot y=x\succ y +x\prec y}
                &&\mbox{pre-Lie algebra}~(A,\ast)\ar[d]^{[x,y]=x\ast y-y\ast x}\\
 \ar[rr]^{\quad [x,y]=x\cdot y-y\cdot x} \mbox{associative algebra}~(A,\cdot)
               && \mbox{Lie algebra}~(A,[-,-]).   }
$$
A Zinbiel algebra can be viewed as a commutative dendriform algebra, namely $x\succ y=y\prec x$. In fact, from the operadic point of view,   dendriform algebras and Zinbiel algebras can be viewed as the splitting of associative algebras and commutative associative algebras respectively (\cite{A2,BBGN,Fard,PBG}).

In this paper, we introduce the notion of a noncommutative pre-Poisson algebra, which consists of a dendriform algebra and a pre-Lie algebra, such that some compatibility conditions are satisfied.  Through the sub-adjacent associative algebra and the sub-adjacent Lie algebra, a noncommutative pre-Poisson algebra gives rise to a noncommutative Poisson algebra naturally. Thus, {noncommutative pre-Poisson algebras can be viewed as the splitting of noncommutative Poisson algebras}.
 We further introduce the notion of a Rota-Baxter operator (more generally an $\huaO$-operator) on a noncommutative Poisson algebra, which is simultaneously a Rota-Baxter operator on the underlying associative algebra and a Rota-Baxter operator on the underlying Lie algebra. See \cite{A2,Bai2007,BGN2013A,CK,EGK,Gub,GK,Rota,Uchino2} for more details on Rota-Baxter operators and $\huaO$-operators. A noncommutative pre-Poisson algebra can be obtained through the action of a Rota-Baxter operator (more generally an $\huaO$-operator).
 The above relations can be summarized into the following
commutative diagram:
$$
\xymatrix{
\mbox{{\bf dendriform algebra} + pre-Lie algebra }\ar[rr] \ar[d]_{\mbox{sub-adjacent~~}}
                && \mbox{noncomm pre-Poisson algebra}\ar[d]_{\mbox{sub-adjacent~~}}\\
\mbox{{\bf associative algebra} + Lie algebra }\ar[rr] \ar@<-1ex>[u]_{\mbox{~~Rota-Baxter ~action}}
                && \mbox{noncomm Poisson algebra}\ar@<-1ex>[u]_{\mbox{~~Rota-Baxter ~action}}}
$$

 We construct skew-symmetric solutions of the Poisson Yang-Baxter  equation  in some special noncommutative Poisson algebras obtained from these structures.

\subsection{Outline of the paper} In Section \ref{sec:quasirep}, we recall quasi-representations and the corresponding cohomology theory of noncommutative Poisson algebras, and introduce the notion of a coherent noncommutative Poisson algebra for our later study of noncommutative Poisson bialgebras.

In Section \ref{bialgebras}, we introduce the notions of matched pairs, Manin triples for noncommutative Poisson algebras and noncommutative (pseudo)-Poisson bialgebras. The equivalences between matched pairs of coherent noncommutative Poisson algebras, Manin triples for noncommutative Poisson algebras and noncommutative Poisson bialgebras are established.

In Section \ref{sec:coboundary}, we study coboundary
noncommutative Poisson bialgebras with the help of
quasi-representations of noncommutative Poisson algebras and the
corresponding cohomology theory, which leads to the introduction
of  the  Poisson Yang-Baxter equation in a coherent
noncommutative Poisson algebra.

In Section \ref{sec:pre-Poisson}, we introduce the notion of  a
noncommutative pre-Poisson algebra and a Rota-Baxter operator
(more generally an $\huaO$-operator) on  a noncommutative Poisson
algebra, by which we construct skew-symmetric solutions of the
Poisson Yang-Baxter  equation in certain
special noncommutative Poisson algebras obtained from these
structures.

In this paper, all the vector spaces are over an algebraically closed field $\mathbb K$ of characteristic $0$, and finite dimensional.

\vspace{2mm}

\section{Quasi-representations and cohomologies of noncommutative Poisson algebras}\label{sec:quasirep}
In this section, we recall (quasi)-representations of noncommutative Poisson algebras and the corresponding cohomology theory.
\begin{defi}
  Let $(A,\cdot_A)$ be an associative algebra and $V$   a vector space. Let $\huaL,\huaR:A\longrightarrow\gl(V)$ be two linear maps with $x\rightarrow \huaL_x$ and $x\rightarrow \huaR_x$ respectively. The triple $(V;\huaL,\huaR)$ is called a {\bf representation} of $A$ if for all  $x,y\in A,$ we have
  $$\huaL_{x\cdot_A y}=\huaL_x\circ \huaL_y,\quad \huaR_{x\cdot_A y}v=\huaR_y\circ \huaR_x,\quad \huaL_x\circ \huaR_y=\huaR_y\circ\huaL_x.$$
 \end{defi}


 In fact, $(V;\huaL,\huaR)$ is a representation of an associative algebra
$A$ if and only if the direct sum $A\oplus V$ of vector spaces is
 an associative algebra (the semi-direct product) by
defining the multiplication on $A\oplus V$ by
$$
  (x_1+v_1)\cdot_{(\huaL,\huaR)}(x_2+v_2)=x_1\cdot_A x_2+\huaL_{x_1}v_2+\huaR_{x_2}v_1,\quad \forall~ x_1,x_2\in A,v_1,v_2\in V.
$$
  We denote it by $A\ltimes_{\huaL,\huaR} V$ or simply by $A\ltimes V$.

\begin{lem}
Let $(V;\huaL,\huaR)$ be a representation of an associative algebra $(A,\cdot_A)$.
   Define $\huaL^*:A\longrightarrow \gl(V^*)$ and $\huaR^*:A\longrightarrow \gl(V^*)$  by
$$
 \langle \huaL^*_x\alpha,v\rangle=-\langle \alpha,\huaL_xv\rangle,\quad\langle \huaR^*_x\alpha,v\rangle=-\langle \alpha,\huaR_xv\rangle,\quad \forall ~ x\in A,\alpha\in V^*,v\in V.
$$
  Then $(V^*;-\huaR^*,-\huaL^*)$ is a representation of $(A,\cdot_A)$.
\end{lem}

\begin{ex}{\rm
  Let $(A,\cdot_A)$ be an associative algebra. Let $L_x$ and $R_x$ denote the left and right multiplication operators, respectively, that is, $L_xy=x\cdot_A y,R_yx=x\cdot_A y$ for all $x,y\in A$. Then   $(A;L,R)$ is a representation of $(A,\cdot_A)$, called the {\bf regular representation}. Furthermore,  $(A^*;-R^*,-L^*)$ is also a representation of $(A,\cdot_A)$.
  }
\end{ex}

 Similarly, let $(\frak g,[-,-]_{\frak g})$ be a Lie algebra  and $V$ a vector space. Let $\rho:\frak g\rightarrow
\frak g\frak l(V)$ be a linear map. The pair $(V;\rho)$ is called a {\bf representation} of $\frak g$ if for all
$x,y\in \frak g$, we have
$$\rho([x,y]_{\frak g})=[\rho(x),\rho(y)].$$
In fact, $(V;\rho)$ is a representation of a Lie algebra
$\frak g$ if and only if the direct sum $\frak g\oplus V$ of vector spaces is
 a Lie algebra (the semi-direct product) by
defining the Lie bracket on $\frak g\oplus V$ by
$$
  [x_1+v_1,x_2+v_2]_\rho=[x_1,x_2]_{\frak g}+\rho(x_1)(v_2)-\rho(x_2)(v_1),\quad \forall~x_1,x_2\in \frak g,v_1,v_2\in V.
$$
  We denote it by $\frak g\ltimes_\rho V$ or simply by $\frak g\ltimes V$. Moreover, let $(V;\rho)$ be a representation of a Lie algebra $(\frak g,[-,-]_{\frak g})$.
   Define $\rho^*:\frak g\longrightarrow \gl(V^*)$  by
$$
 \langle \rho^*(x)(\alpha),v\rangle=-\langle \alpha,\rho(x) (v)\rangle,\quad \forall~x\in \frak g,\alpha\in V^*,v\in V.
$$
  Then $(V^*;\rho^*)$ is a representation of $(\g, [-,-]_\g)$. In particular, define $\ad:\g\lon\gl(\g)$ by $\ad_xy=[x,y]_\g$ for all $x,y\in \g$. Then   $(\frak g;\ad)$ is a representation of $(\g, [-,-]_\g)$, called the {\bf adjoint representation}. Furthermore,  $(\g^*;\ad^*)$ is also a representation of $(\g, [-,-]_\g)$.

\begin{defi}{\rm{(\cite{FGV,YYY})}} Let $(P,\cdot_P,\{-,-\}_P)$ be a noncommutative Poisson algebra.
\begin{itemize}
 \item[{\rm(i)}]   A {\bf quasi-representation} of $P$ is a quadruple $(V;\huaL,\huaR,\rho)$ such that $(V;\huaL,\huaR)$ is a representation of the associative algebra $(P,\cdot_P)$ and $(V;\rho)$ is a  representation of the Lie algebra $(P,\{-,-\}_P)$ satisfying
   \begin{eqnarray}
     \label{eq:rep 1}\huaL_{{\{x,y\}}_P}&=&\rho(x)\huaL_y-\huaL_y\rho(x),\\
     \label{eq:rep 2} \huaR_{\{x,y\}_P}&=&\rho(x)\huaR_y-\huaR_y\rho(x),\quad \forall~x, y\in P.
   \end{eqnarray}
 \item[{\rm(ii)}] A quasi-representation $(V;\huaL,\huaR,\rho)$ is called a {\bf representation} of $P$  if we have
   \begin{equation}\label{eq:con-rep}
     \rho(x\cdot_P y)=\huaL_x\rho(y)+\huaR_y\rho(x).
   \end{equation}
   \end{itemize}
  \end{defi}

 By a direct calculation, we have
\begin{pro}
   Let $(P,\cdot_P,\{-,-\}_P)$ be a noncommutative Poisson algebra.
   \begin{itemize}
 \item[{\rm(i)}] If $(V;\huaL,\huaR,\rho)$ is a quasi-representation of  $P$, then $(V^*;-\huaR^*,-\huaL^*,\rho^*)$ is also a quasi-representation of  $P$;

 \item[{\rm(ii)}]  If $(V;\huaL,\huaR,\rho)$ is a quasi-representation and satisfies
    \begin{equation}\label{eq:coherent condition}
       \rho(x\cdot_P y)=\rho(x)\huaL_y+\rho(y)\huaR_x,
   \end{equation}
   then $(V^*;-\huaR^*,-\huaL^*,\rho^*)$ is a representation of $P$.
   \end{itemize}
\end{pro}

\begin{ex}{\rm
   Let $(P,\cdot_P,\{-,-\}_P)$ be a noncommutative Poisson algebra. Then $(P;L,R,\ad)$ is a representation of $P$, which is also called the {\bf regular representation} of  $P$. However, the dual  $(P^*;-R^*,-L^*,\ad^*)$ is just a quasi-representation of $P$. Furthermore, it is straightforward to check that $(P^*;-R^*,-L^*,\ad^*)$ is a representation of $P$ if and only if the noncommutative Poisson algebra $P$ satisfies
   \begin{equation}\label{eq:coboundary extra condition}
     \Poisson{x,y\cdot_P z}_P+\Poisson{y,z\cdot_P x}_P+ \Poisson{z,x\cdot_P y}_P=0,\quad\forall~x,y\in P.
  \end{equation}
  }
\end{ex}

\begin{defi}
A  noncommutative Poisson algebra $(P,\cdot_P,\{-,-\}_P)$ is called {\bf coherent} if it satisfies \eqref{eq:coboundary extra condition}.
\end{defi}

\begin{pro}
 A noncommutative Poisson algebra $(P,\cdot_P,\{-,-\}_P)$ is  coherent if and only if it satisfies
 \begin{eqnarray*}
   [\{x,y\}_P,z]_P+ [\{z,x\}_P,y]_P+ [\{y,z\}_P,x]_P=0
 \end{eqnarray*}
 for all $x,y,z\in P$, where $[-,-]_P:P\times P\longrightarrow P$ is the commutator Lie bracket defined by
 \begin{equation}
 \label{eq:commutator bracket}[x,y]_P=x\cdot_P y-y\cdot_P x.
 \end{equation}
\end{pro}
\pf By the Leibniz rule of the noncommutative Poisson algebra $(P,\cdot_P,\{-,-\}_P)$, we have
\begin{eqnarray*}
   &&\Poisson{x,y\cdot_P z}_P+\Poisson{y,z\cdot_P x}_P+\Poisson{z,x\cdot_P y}_P\\
   &=&\Poisson{x,y}_P\cdot_P z+y\cdot_P\Poisson{x,z}_P+\Poisson{y,z}_P\cdot_P x+z\cdot_P\Poisson{y,x}_P+\Poisson{z,x}_P\cdot_P y+x\cdot_P\Poisson{z,y}_P\\
   &=&\Poisson{x,y}_P\cdot_P z-z\cdot_P\Poisson{x,y}_P+\Poisson{z,x}_P\cdot_P y-y\cdot_P\Poisson{z,x}_P+\Poisson{y,z}_P\cdot_P x-x\cdot_P\Poisson{y,z}_P\\
   &=& [\{x,y\}_P,z]_P+ [\{z,x\}_P,y]_P+ [\{y,z\}_P,x]_P.
\end{eqnarray*}
Then the conclusion follows immediately.\qed\vspace{3mm}

  Similarly,  a coherent noncommutative Poisson
algebra $(P,\cdot_P,\{-,-\}_P)$ also satisfies
 $$\{[x,y]_P,z\}_P+ \{[z,x]_P,y\}_P+ \{[y,z]_P,x\}_P=0.$$

\begin{cor}
Any commutative Poisson algebra is coherent.
\end{cor}

 Coherent noncommutative Poisson algebras are closely related to compatible Lie algebras (\cite{GS1,GS2,GS3,OS}).

\begin{defi}{ A {\bf compatible Lie algebra} $(\g,[-,-]_{1},[-,-]_{2})$
consists of two Lie algebras $(\g,[-,-]_{1})$  and $(\g,[-,-]_{2})$ such that for any $k_1,k_2\in\mathbb{K}$, the following bilinear operation
\begin{equation}
\label{eq:Lie} [x,y]=k_1[x,y]_{1}+k_2[x,y]_{2},\quad \forall~x,y\in\mathfrak{g}
\end{equation}
defines a Lie algebra structure on $\g$.
}\end{defi}

\begin{pro} {\rm (\cite{GS1}) }
 Let $(\g,[-,-]_{1})$ and $(\g,[-,-]_{2})$ be two Lie algebras.
 Then  $(\g,[-,-]_{1},[-,-]_{2})$ is a compatible Lie algebra
if and only if for any $x,y,z\in \mathfrak{g}$, the following equation holds
\begin{equation}
[[x,y]_{1},z]_{2}+[[y,z]_{1},x]_{2}+[[z,x]_{1},y]_{2}
+[[x,y]_{2},z]_{1}+[[y,z]_{2},x]_{1}+[[z,x]_{2},y]_{1}=0.\label{eq:cl}
\end{equation}
\end{pro}

\begin{cor}
Let $(P,\cdot_P,\{-,-\}_P)$ be a coherent noncommutative Poisson algebra. Then we have a compatible Lie algebra $(P, \{-,-\}_P,[-,-]_P)$, where the bracket $[-,-]_P$ is given by \eqref{eq:commutator bracket}.
\end{cor}

\begin{ex}\label{ex:coherent poisson 2}{\rm
  Let $(A,\cdot_A)$ be an associative algebra. Define a bracket $\Poisson{-,-}_\hbar$ as follows
  \begin{equation}
    \Poisson{x,y}_\hbar=\hbar(x\cdot_A y-y\cdot_A x),
  \end{equation}
  where $x,y\in A$ and $\hbar$ is a fixed number. Then $(A,\cdot_A,\Poisson{-,-}_\hbar)$ is a coherent noncommutative Poisson algebra, which is called the {\bf standard noncommutative Poisson algebra}.
  }
\end{ex}

By a direct calculation, we have
\begin{pro}
  The standard noncommutative Poisson algebra $(P,\cdot_P,\Poisson{-,-}_\hbar)$ is compatible with any coherent noncommutative Poisson algebra $(P,\cdot_P,\{-,-\}_P)$   in the sense of that for any $k_1,k_2\in \K$,  $(P,k_1\cdot_P+k_2\cdot_P,k_1\{-,-\}_P+k_2\Poisson{-,-}_\hbar)$ is still a coherent noncommutative Poisson algebra.
\end{pro}

\begin{ex}\label{ex:SNCPA1}{\rm
  Let $P$ be a $3$-dimensional vector space  with basis $\{e_1,e_2,e_3\}$. Define the non-zero multiplication and the bracket operation  on $P$ by
  \begin{eqnarray*}
  e_1\cdot e_2&=&e_3,\quad e_2\cdot e_1=-e_3;\\
 \Poisson{e_1,e_2}&=&a e_1+b e_2+c e_3,\quad \Poisson{e_1,e_3}=b e_3,\quad\Poisson{e_2,e_3}=-a e_3,\quad\forall~a,b,c \in \K.
  \end{eqnarray*}
  Then $(P,\cdot,\Poisson{-,-})$ is a coherent noncommutative Poisson algebra.
  }
\end{ex}

\begin{ex}\label{ex:SNCPA2}{\rm
  Let $P$ be a $4$-dimensional vector space with basis $\{e_1,e_2,e_3,e_4\}$. Define the non-zero multiplication and the   non-zero bracket  operation on $P$ by
  \begin{eqnarray*}
  e_1\cdot e_1&=&e_4,\quad e_1\cdot e_2=e_4,\quad e_3\cdot e_2=e_4,\quad e_3\cdot e_3=-e_4;\\
 \Poisson{e_1,e_2}&=&a e_4,\quad \Poisson{e_1,e_3}=b e_4,\quad\Poisson{e_2,e_3}=c e_4,\quad\forall~a,b,c \in \K.
  \end{eqnarray*}
  Then $(P,\cdot,\Poisson{-,-})$ is a coherent noncommutative Poisson algebra.
  }
\end{ex}

It is straightforward to obtain the following conclusion.
\begin{pro}\label{pro:semi-direct}
 Let $(P,\cdot_P,\{-,-\}_P)$ be a noncommutative Poisson algebra and $(V;\huaL,\huaR,\rho)$ a representation. Then $(P\oplus V,\cdot_{(\huaL,\huaR)},\Poisson{-,-}_\rho)$ is a noncommutative Poisson algebra, where $(P\oplus V,\cdot_{(\huaL,\huaR)})$ is the semi-direct product associative algebra $P\ltimes_{(\huaL,\huaR)} V$ and $(P\oplus V,\Poisson{-,-}_\rho)$ is the semi-direct product Lie algebra $P\ltimes_{\rho} V$.

 Furthermore, if the  noncommutative Poisson algebra $P$ is coherent and the representation
  $(V;\huaL,\huaR,\rho)$ satisfies \eqref{eq:coherent condition}, then $(P\oplus V,\cdot_{(\huaL,\huaR)},\Poisson{-,-}_\rho)$ is a coherent noncommutative Poisson algebra.
\end{pro}

We denote a semi-direct product noncommutative Poisson algebra by $P\ltimes_{(\huaL,\huaR,\rho)} V$.

\emptycomment{\begin{defi}
We call the semi-direct product $(P\oplus V,\cdot_{(\huaL,\huaR)},\Poisson{-,-}_\rho)$ structure a {\bf quasi-semi-direct product Poisson algebra} if  $(\huaL,\huaR,\rho)$ is quasi-representation of Poisson algebra $(P,\cdot_P,\{-,-\}_P)$ on $V$.
\end{defi}}

The cohomology complex for a noncommutative Poisson algebra $(P,\cdot_P,\{-,-\}_P)$  associated to a quasi-representation $(V;\huaL,\huaR,\rho)$ is given as follows (\cite{CohomologyPA1}). Denote
$$C^{i,j}(P,V)=\Hom( (\otimes^iP)\otimes(\wedge^jP), V)\quad\mbox{and}\quad C^n(P,V)=\sum_{i+j=n}C^{i,j}(P,V).$$
 The coboundary operator $\delta^n:C^{n}(P,V)\longrightarrow C^{n+1}(P,V)$ is given by
 $$\delta^n=\sum_{i+j=n}(\hat{\delta}^{i,j}+(-1)^i\bar{\delta}^{i,j}),$$
 where $\hat{\delta}^{i,j}:C^{i,j}(P,V)\longrightarrow C^{i+1,j}(P,V)$ is given by
 \begin{eqnarray*}
  &&\hat{\delta}^{i,j}\varphi\big((a_1\otimes\cdots \otimes a_{i+1})\otimes(x_1\wedge\cdots\wedge x_j)\big)\\
  &=&\huaL_{a_1} \varphi((a_2\otimes\cdots \otimes a_{i+1})\otimes(x_1\wedge\cdots\wedge x_j))\\
  &&+\sum_{k=1}^i(-1)^k\varphi\big((a_1\otimes\cdots\otimes a_k\cdot_P a_{k+1}\otimes\cdots \otimes a_{i+1})\otimes(x_1\wedge\cdots\wedge x_j)\big)\\
&&+(-1)^{i+1}\huaR_{ a_{i+1}}\varphi\big((a_1\otimes\cdots \otimes a_{i})\otimes(x_1\wedge\cdots\wedge x_j)\big)
 \end{eqnarray*}
 and
$\bar{\delta}^{i,j}:C^{i,j}(P,V)\longrightarrow C^{i,j+1}(P,V)$ is given by
 \begin{eqnarray*}
  &&\bar{\delta}^{i,j}\varphi\big((a_1\otimes\cdots \otimes a_{i})\otimes(x_1\wedge\cdots\wedge x_{j+1})\big)\\
  &=&\sum_{l=1}^{j+1}(-1)^{l+1}\Big(\rho(x_l)\varphi\big((a_1\otimes\cdots \otimes a_{i})\otimes(x_1\wedge\cdots \hat{x}_l\cdots\wedge x_{j+1})\big)\\
 &&-\sum_{k=1}^i\varphi\big((a_1\otimes\cdots\otimes \Poisson{x_l,a_k}_P\otimes\cdots \otimes a_{i})\otimes(x_1\wedge\cdots \hat{x}_l\cdots\wedge x_{j+1})\big)\Big)\\
 &&+\sum_{1\leq p\leq q\leq j+1}(-1)^{p+q}\varphi\big((a_1\otimes\cdots \otimes a_{i})\otimes(\Poisson{x_p,x_q}_P\wedge x_1\wedge \cdots \hat{x}_p\cdots \hat{x}_q\cdots\wedge x_{j+1})\big)
 \end{eqnarray*}
 for all $\varphi\in C^{i,j}(P,V) $ and $a_1,\cdots,a_{i+1}, x_1,\cdots,x_{j+1}\in P$.

In particular, a 1-cochain $\varphi+\psi$, where $\varphi\in
C^{1,0}(P,V)$ and $\psi\in C^{0,1}(P,V)$, is a $1$-cocycle on $P$
with   coefficients in the quasi-representation
$(V;\huaL,\huaR,\rho)$ means that $\delta^1(\varphi+\psi)=0$, i.e.
 $$\hat{\delta}^{1,0} \varphi(a\otimes b)=0,\quad\big(\hat{\delta}^{0,1} \psi-\bar{\delta}^{1,0}\varphi\big)(a\otimes x)=0, \quad\bar{\delta}^{0,1} \psi(x\otimes y)=0,\quad \forall~a,b,x,y\in P.$$
More precisely,
\begin{eqnarray}
  \varphi(a\cdot_P b)&=&\huaL_{a}\varphi(b)+\huaR_{b}\varphi(a);\\
  \varphi(\Poisson{x,a}_P)&=&\rho(x)  \varphi(a)-\huaL_{a}\psi(x)+\huaR_{a}\psi(x);\\
  \psi(\Poisson{x,y}_P)&=&\rho(x) \psi(y)-\rho(y)  \psi(x).
\end{eqnarray}
A 1-cochain $\varphi+\psi\in C^1(P,V)$ is a $1$-coboundary on $P$ with   coefficients in the  quasi-representation $(V;\huaL,\huaR,\rho)$ if and only if there exists an element $u\in V$ such that
\begin{equation}
  \huaL_{a}u-\huaR_{a}u=\varphi(a)\quad\mbox{and}\quad \rho(x)u=\psi(x),\quad\forall~a,x\in P.
\end{equation}

\emptycomment{Furthermore, a $2$-cochain $\varphi+\phi+\psi$, where $\varphi\in C^{2,0}(P,V)$, $\phi\in C^{1,1}(P,V)$ and $\psi\in C^{0,2}(P,V)$, is a 2-cocycle on  $P$ with   coefficients in the quasi-representation $(V;\huaL,\huaR,\rho)$ means that
  $\delta^2(\varphi+\phi+\psi)=0$, i.e. {\bf CM: i.e.,} for $a_1,a_2,a_3,x_1,x_2,x_3\in P$, there holds:
\begin{eqnarray*}
\hat{\delta}^{2,0} \varphi(a_1\otimes a_2\otimes a_3)=0,&&\big(\hat{\delta}^{1,1} \phi+\bar{\delta}^{2,0}\varphi\big)(a_1\otimes a_2\otimes x_1)=0, \\
\bar{\delta}^{0,2} \psi(x_1\wedge x_2\wedge x_3)=0,&& (\hat{\delta}^{0,2}\psi-\bar{\delta}^{1,1} \phi)(a_1\otimes x_1\wedge x_2)=0.
\end{eqnarray*}}

\emptycomment{Let $(A,\succ,\prec)$ be a dendriform algebra. For $x\in A$, define $L^\succ_x,L^\prec_x,R^\prec_x:A\longrightarrow\gl(A)$ by
$$L^\succ_x(y)=x\succ y,\quad L^\prec_x(y)=x\prec y\quad R^\prec_x(y)=y\prec x,\quad\forall~y\in A.$$
Then $(A;L^\succ,R^\prec)$ is a representation of the associative algebra $(A,\cdot)$, where $x\cdot y=x\succ y+x\prec y$. and $(A;L^\succ,-L^\prec)$ is a representation of the adjacent pre-Lie algebra $(A,\ast)$.}

\section{Matched pairs, Manin triples   and   (pseudo-)Poisson bialgebras}\label{bialgebras}
A {\bf matched pair of Lie algebras} is a pair of Lie algebras $(\g_1,[-,-]_{\g_1})$ and $(\g_2,[-,-]_{\g_2})$ together with two representations $\rho:\g_1\longrightarrow\gl(\g_2)$ and $\varrho:\g_2\longrightarrow\gl(\g_1)$ satisfying
\begin{eqnarray*}
 \varrho(\alpha){[x,y]_{\g_1}}&=&[ \varrho(\alpha)x,y]_{\g_1}+[x, \varrho(\alpha)y]_{\g_1}- \varrho(\rho(x)\alpha)y+ \varrho(\rho(y)\alpha)x,\\
   \rho(x)[\alpha,\beta]_{\g_2}&=&[\rho(x)\alpha,\beta]_{\g_2}+[\alpha,\rho(x)\beta]_{\g_2}-\rho(\varrho(\alpha)x)\beta+\rho(\varrho(\beta)x)\alpha
\end{eqnarray*}
for all $x,y\in\g_1$ and $\alpha,\beta\in\g_2$. In this case, there exists a Lie algebra structure on the vector space $\frkd=\g_1\oplus\g_2$ given by
\begin{equation}
  [x+\alpha,y+\beta]_\frkd=[x,y]_{\g_1}+\varrho(\alpha)y-\varrho(\beta)x+[\alpha,\beta]_{\g_2}+\rho(x)\beta-\rho(y)\alpha.
\end{equation}
It is denoted by $\g_1\bowtie^{\varrho}_{\rho}\g_2$ or simply by $\g_1\bowtie\g_2$, and called the double of the matched pair. Moreover, every Lie algebra which is a direct sum of the underlying vector spaces of two subalgebras is the double of a matched pair of Lie algebras.

A {\bf matched pair of associative algebras} is a pair of associative algebras $(A_1,\cdot_1)$ and $(A_2,\cdot_2)$ together with two representations $(\huaL,\huaR):A_1\longrightarrow\gl(A_2)$ and $(\tilde{\huaL},\tilde{\huaR}):A_2\longrightarrow\gl(A_1)$ satisfying
\begin{eqnarray*}
&& \huaL_x(\alpha\cdot_2\beta)=\huaL_{\tilde{\huaR}_\alpha x}\beta+(\huaL_x\alpha)\cdot_2\beta,\qquad \huaR_x(\alpha\cdot_2\beta)=\huaR_{\tilde{\huaL}_\beta x}\alpha+\alpha\cdot_2\huaR_x\beta,\\
&&\tilde{\huaL}_\alpha(x\cdot_1y)=\tilde{\huaL}_{\huaR_x \alpha}y+(\tilde{\huaL}_\alpha x)\cdot_1y,\qquad \tilde{\huaR}_\alpha(x\cdot_1y)=\tilde{\huaR}_{\huaL_y \alpha}x+ x\cdot_1\tilde{\huaR}_\alpha y,
\\
&&\huaL_{\tilde{\huaL}_\alpha x}\beta+(\huaR_x\alpha)\cdot_2 \beta-\huaR_{\tilde{\huaR}_\beta x}\alpha-\alpha\cdot_2 \huaL_x\beta=0,\quad
\tilde{\huaL}_{\huaL_x \alpha}y+(\tilde{\huaR}_\alpha x)\cdot_1 y-\tilde{\huaR}_{\huaR_y \alpha}x-x\cdot_1 \tilde{\huaL}_\alpha y=0
\end{eqnarray*}
for all $x,y\in A_1$ and $\alpha,\beta\in A_2$. In this case, there exists an associative algebra structure on the vector space $\huaA=A_1\oplus A_2$ given by
\begin{equation}
  (x+\alpha)\cdot_\huaA (y+\beta)=x\cdot_1 y+\tilde{\huaL}_\alpha y+\tilde{\huaR}_\beta x+\alpha\cdot_2\beta+\huaL_x \beta+\huaR_y\alpha.
\end{equation}
It is denoted by $A_1\bowtie^{(\tilde{\huaL},\tilde{\huaR})}_{(\huaL,\huaR)}A_2$ or simply by $A_1\bowtie A_2$, and called the double of the matched pair. Moreover, every associative algebra which is a direct sum of the underlying vector spaces of two subalgebras is the double of a matched pair of associative algebras.

\begin{defi}
Let $(P_1,\cdot_1,\{-,-\}_{P_1})$ and $(P_2,\cdot_2,\{-,-\}_{P_2})$ be two noncommutative Poisson algebras. A {\bf matched pair of noncommutative Poisson algebras} is a pair of noncommutative Poisson algebras $(P_1,P_2)$ together with two representations $(\huaL,\huaR,\rho):P_1\longrightarrow\gl(P_2)$ and $(\tilde{\huaL},\tilde{\huaR},\varrho):P_2\longrightarrow\gl(P_1)$ such that $P_1\bowtie^{\varrho}_{\rho}P_2$ is a matched pair of Lie algebras, $P_1\bowtie^{(\tilde{\huaL},\tilde{\huaR})}_{(\huaL,\huaR)}P_2$ is a matched pair of associative algebras and for all $x,y\in P_1$ and $\alpha,\beta\in P_2$, the following equalities hold:
\begin{eqnarray}
\label{eq:MPP1}  \rho(x)(\alpha\cdot_2\beta)&=&(\rho(x)\alpha)\cdot_2 \beta+\alpha\cdot_2 \rho(x)\beta-\huaL_{\varrho(\alpha)x}\beta-\huaR_{\varrho(\beta)x}\alpha,\\
 \label{eq:MPP2} \huaL_x\Poisson{\alpha,\beta}_{P_2}&=&\Poisson{\alpha,\huaL_x\beta}_{P_2}-\rho(\tilde{\huaR}_\beta x)\alpha-\huaL_{\varrho(\alpha)x}\beta+(\rho(x)\alpha)\cdot_2 \beta,\\
\label{eq:MPP3}  \varrho(\alpha)(x\cdot_1y)&=&(\varrho(\alpha)x)\cdot_1 y+x\cdot_1 \varrho(\alpha)y-\tilde{\huaL}_{\rho(x)\alpha}y-\tilde{\huaR}_{\rho(y)\alpha}x,\\
\label{eq:MPP4} \tilde{ \huaL}_\alpha\Poisson{x,y}_{P_1}&=&\Poisson{x,\tilde{\huaL}_\alpha y}_{P_1}-\varrho(\huaR_y \alpha)x-\tilde{\huaL}_{\rho(x)\alpha}y+(\varrho(\alpha)x)\cdot_1 y.
\end{eqnarray}

\end{defi}

\emptycomment{\begin{defi}
A {\bf quasi-noncommutative Poisson algebra} is a triple $(P,\cdot_P,\{-,-\}_P)$, where $(P,\cdot_P)$ is an associative algebra, $(P,\{-,-\}_P)$ is a Lie algebra and there exists a decomposition of $P=P_1\oplus P_2$ as vector spaces such that $(P_1,\cdot_P,\Poisson{-,-}_{P})$ and $(P_2,\cdot_P,\Poisson{-,-}_{P})$ are noncommutative Poisson algebras and satisfies
  \begin{eqnarray}
    \{x,y\cdot_P u\}_P&=&\{x,y\}_P\cdot_P u+y\cdot_P\{x,u\}_P,\\
    \{x,u\cdot_P y\}_P&=&\{x,u\}_P\cdot_P y+u\cdot_P\{x,y\}_P,\\
    \{u,v\cdot_P x\}_P&=&\{u,v\}_P\cdot_P x+v\cdot_P\{u,x\}_P,\\
    \{u,x\cdot_P v\}_P&=&\{u,x\}_P\cdot_P v+x\cdot_P\{u,v\}_P
  \end{eqnarray}
for all $x,y\in P_1$ and $u,v\in P_2$. We call $(P_1,P_2)$ the noncommutative Poisson algebra pair associated to $P$.
\end{defi}}

\begin{pro}\label{pro:Matched pair and double}
  Let $(P_1,\cdot_1,\{-,-\}_{P_1})$ and $(P_2,\cdot_2,\{-,-\}_{P_2})$ be two noncommutative Poisson algebras. If $(P_1,P_2)$ is a matched pair of noncommutative Poisson algebras with representations $(\huaL,\huaR,\rho):P_1\longrightarrow\gl(P_2)$ and $(\tilde{\huaL},\tilde{\huaR},\varrho):P_2\longrightarrow\gl(P_1)$, then there exists a noncommutative Poisson algebra structure on $P=P_1\oplus P_2$ given by
  \begin{eqnarray}
    \Poisson{x+\alpha,y+\beta}_P&=&\Poisson{x,y}_{P_1}+\varrho(\alpha)y-\varrho(\beta)x+\Poisson{\alpha,\beta}_{P_2}+\rho(x)\beta-\rho(y)\alpha,\\
    (x+\alpha)\cdot_P (y+\beta)&=&x\cdot_1 y+\tilde{\huaL}_\alpha y+\tilde{\huaR}_\beta x+\alpha\cdot_2\beta+\huaL_x \beta+\huaR_y\alpha
  \end{eqnarray}
  for all $x,y\in P_1$ and $\alpha,\beta\in P_2$.
\end{pro}
\pf The verification is a routine calculation and thus omitted.\qed\vspace{3mm}

We denote this noncommutative Poisson algebra by $P_1\bowtie^{(\tilde{\huaL},\tilde{\huaR},\varrho)}_{(\huaL,\huaR,\rho)}P_2$ or simply by $P_1\bowtie P_2$, and call it {\bf the double} of the matched pair. Moreover, every noncommutative Poisson algebra which is a direct sum of the underlying vector spaces of two  subalgebras is the double of a matched pair of noncommutative Poisson algebras.

A {\bf quadratic  noncommutative Poisson algebra} is a noncommutative Poisson algebra $P$ equipped with a nondegenerate symmetric bilinear form $\frkB(-,-)$ on $P$ which is {\bf invariant} in the sense that
  \begin{equation}
    \frkB(\Poisson{a,b}_P,c)=\frkB(a,\Poisson{b,c}_P),\quad \frkB(a\cdot_P b,c)=\frkB(a,b\cdot_P c),\quad\forall~a,b,c\in P.
  \end{equation}

\begin{defi}\label{defi:Manin}
  A {\bf  Manin triple of noncommutative Poisson algebras} is a triple $((P,\frkB),P_1,P_2)$, where $(P,\huaB)$ is a quadratic noncommutative Poisson algebra, $(P_1,\cdot_1,\Poisson{-,-}_{P_1})$ and $(P_2,\cdot_2,\Poisson{-,-}_{P_2})$ are noncommutative Poisson subalgebras of $P$,
  such that
  \begin{enumerate}
\item[\rm(i)]  $P=P_1\oplus P_2$ as vector spaces;
\item[\rm(ii)] $P_1$ and $P_2$ are isotropic with respect to $\frkB(-,-)$.
\end{enumerate}
 \end{defi}
\begin{rmk}\label{rmk:Manin}
 It is obvious that a Manin triple of noncommutative Poisson algebras is simultaneously a Manin triple of Lie algebras (\cite{CP}) and an analogue of a Manin triple in the context of associative algebras which is called
  a double construction of Frobenius algebras in \cite{Bai1}.
\end{rmk}

\begin{pro}\label{thm:Manin-Standard}
Let $(P,P_1,P_2)$ be a Manin triple of noncommutative Poisson algebras. Then the noncommutative Poisson algebras $P$, $P_1$ and $P_2$ must be coherent.
\end{pro}
\pf First we show that the noncommutative Poisson algebra $P_1$ is coherent. Since $P_1$ is a subalgebra of $P$, we only need to prove that for any $x,y,z\in P_1$, the following equality holds:
\begin{eqnarray}\label{eq:Manin eq1}
  \Poisson{x\cdot_P y,z}_P=\Poisson{x,y\cdot_P z}_P+\Poisson{y,z\cdot_P x}_P.
\end{eqnarray}
By the invariance of $\frkB$ and the Leibniz rule, for $\alpha\in P_2$, we have
\begin{eqnarray*}
  \frkB\big(\Poisson{x\cdot_P y,z}_P-\Poisson{x,y\cdot_P z}_P-\Poisson{y,z\cdot_P x}_P,\alpha\big)
  =\frkB\big(y \cdot_P\Poisson{z,\alpha}_P-\Poisson{y\cdot_P z,\alpha}_P+\Poisson{y,\alpha}_P\cdot_P z,x\big)=0.
\end{eqnarray*}
  Then by the nondegeneracy of $\frkB$, \eqref{eq:Manin eq1} follows immediately.

Also by the invariance of $\frkB$, for all  $x,y,z\in P_1$ and $\alpha,\beta\in P_2$,  we have
 \begin{eqnarray*}
  \frkB\big(\Poisson{x\cdot_P y,\alpha}_P-\Poisson{x,y\cdot_P \alpha}_P-\Poisson{y,\alpha\cdot_P x}_P,z\big)
  =\frkB\big(\Poisson{z,x\cdot_P y}_P-\Poisson{z,x}_P\cdot_P y-x\cdot_P\Poisson{z,y}_P,\alpha\big)=0,
\end{eqnarray*}
 and
\begin{eqnarray*}
 \frkB\big(\Poisson{x\cdot_P y,\alpha}_P-\Poisson{x,y\cdot_P \alpha}_P-\Poisson{y,\alpha\cdot_P x}_P,\beta\big)
&=&-\frkB\big(\Poisson{x\cdot_P y,\beta}_P-\Poisson{x,\beta}_P\cdot_P y-x\cdot_P\Poisson{y,\beta}_P,\alpha\big)\\
&=&0,
\end{eqnarray*}
which imply that
$$\Poisson{x\cdot_P y,\alpha}_P=\Poisson{x,y\cdot_P \alpha}_P+\Poisson{y,\alpha\cdot_P x}_P.$$

Similarly, for $x\in P_1$ and $\alpha,\beta,\gamma\in P_2$, we also have
\begin{eqnarray*}
  \Poisson{\alpha\cdot_P \beta,\gamma}_P&=&\Poisson{\alpha,\beta\cdot_P \gamma}_P+\Poisson{\beta,\gamma\cdot_P \alpha}_P,\\
  \Poisson{\alpha\cdot_P \beta,x}_P&=&\Poisson{\alpha,\beta\cdot_P x}_P+\Poisson{\beta,x\cdot_P \alpha}_P.
\end{eqnarray*}
Thus the noncommutative Poisson algebras $P$ and $P_2$ are coherent.  \qed\vspace{2mm}

We recall the notions of Lie bialgebras and antisymmetric infinitesimal bialgebras before we give the notion of a noncommutative Poisson bialgebra.

A {\bf Lie bialgebra} is a pair $((\g,[-,-]_\g),\delta)$, where $(\g,[-,-]_\g)$ is a Lie algebra and $\delta:\g\longrightarrow \wedge^2\g$ is a linear map such that $(\g,\delta)$ is a Lie coalgebra and $\delta$ is a $1$-cocycle on $\g$ with   coefficients in the representation $(\g\otimes\g;\ad\otimes 1+1\otimes\ad)$, i.e.
\begin{equation}
  \delta([x,y]_\g)=(\ad_x\otimes 1+1\otimes\ad_x)\delta(y)-(\ad_y\otimes 1+1\otimes\ad_y)\delta(x),\quad\forall~x,y\in\g.
\end{equation}
In particular, a Lie bialgebra $(\g,\delta)$ is called {\bf coboundary} if $\delta$ is a $1$-coboundary, that is, there exists an $r\in\g\otimes \g$ such that \begin{equation}\label{eq:coboundary Lie cobracket}
\delta(x)=(\ad_x\otimes 1+1\otimes \ad_x) r,\quad\forall~x,y\in\g.
\end{equation}
A coboundary Lie bialgebra is usually denoted by $(\g,r)$.

An {\bf antisymmetric infinitesimal bialgebra} is a pair $((A,\cdot_A),\Delta)$, where $(A,\cdot_A)$ is an associative algebra and $\Delta:A\longrightarrow A\otimes A$ is a linear map such that $(A,\Delta)$ is a coassociative coalgebra and $\Delta$ is a $1$-cocycle on $A$ with  coefficients in the representation $(A \otimes A;1\otimes L,R\otimes 1)$, i.e.
\begin{eqnarray}
 \Delta(x\cdot_A y)&=&(1\otimes L_x)\Delta(y)+( R_y\otimes 1)\Delta(x),\quad\forall~x,y\in A
\end{eqnarray}
satisfying
\begin{equation}
(L_y\otimes 1-1\otimes R_y)\Delta(x)+\tau\Big((L_x\otimes 1-1\otimes R_x)\Delta(y)\Big)=0,
\end{equation}
where $\tau:A\otimes A\longrightarrow A\otimes A$ is the exchange operator defined by
 \begin{equation}\label{eq:exchange operator}
 \tau(x\otimes y)=y\otimes x,\quad\forall~x,y\in A.
 \end{equation}
   In particular, an  antisymmetric infinitesimal bialgebra $(A,\Delta)$ is called coboundary if $\Delta$ is a $1$-coboundary, that is, there exists an $r\in A\otimes A$ such that \begin{equation}\label{eq:coboundary associative cobracket}
\Delta(x)=(1\otimes L_x -R_x\otimes 1 )r,\quad\forall~x,y\in A.
\end{equation}
A coboundary  antisymmetric infinitesimal bialgebra is usually denoted by $(A,r)$.

\begin{defi}\label{defi:quasiPBA}
A {\bf noncommutative pseudo-Poisson bialgebra} is a triple $((P,\cdot_P,\Poisson{-,-}_P),\Delta,\delta)$, where $(P,\cdot_P,\Poisson{-,-}_P)$ is a coherent noncommutative Poisson algebra, $\Delta:P\longrightarrow P\otimes P$ and $\delta:P\longrightarrow \wedge^2P$ are linear maps such that
  \begin{enumerate}
\item[\rm(a)] $(P^*,\Delta^*,\delta^*)$ is a noncommutative Poisson algebra, where  $\delta^*:\wedge^2 P^*\longrightarrow P^*$ and $\Delta^*:P^*\otimes P^*\longrightarrow P^*$ are the dual maps of $\Delta$ and $\delta$, defined by
    \begin{eqnarray}
     \label{eq:PBA dualopration}\langle\Delta^*(\alpha,\beta),x\rangle=\langle\Delta(x),\alpha\otimes \beta \rangle ,\quad \langle\delta^*(\alpha,\beta),x\rangle=\langle\delta(x),\alpha\wedge \beta \rangle,\quad\forall~x\in P,\alpha,\beta\in P^*;
    \end{eqnarray}

\item[\rm(b)]   $\delta+\Delta\in C^{0,1}(P,P\otimes P)\oplus
C^{1,0}(P,P\otimes P)$ is a $1$-cocycle on $P$ associated to the
quasi-representation $(P\otimes P;1\otimes L,R\otimes 1,\ad\otimes
1+1\otimes \ad)$, i.e.
    \begin{eqnarray}
  \label{eq:PBA 1}\Delta(x\cdot_P y)&=&(1\otimes L_x)\Delta(y)+(R_y\otimes 1)\Delta(x),\\
  \label{eq:PBA 2}\Delta(\Poisson{x,y}_P)&=&(\ad_x\otimes 1+1\otimes \ad_x) \Delta(y)-(1\otimes L_y)\delta(x)+(R_y\otimes 1)\delta(x),\\
 \label{eq:PBA 3}\delta(\Poisson{x,y}_P)&=&(\ad_x\otimes 1+1\otimes \ad_x)\delta(y)-(\ad_y\otimes 1+1\otimes \ad_y)\delta(x);
\end{eqnarray}
\item[\rm(c)] $\triangle$ and $\delta$ satisfy
\begin{eqnarray}
 \label{eq:PBA 4} (L_y\otimes 1-1\otimes R_y)\Delta(x)+\tau\Big((L_x\otimes 1-1\otimes R_x)\Delta(y)\Big)&=&0,\\
 \label{eq:PBA 5}(L_x\otimes 1)\delta(y)+(R_y\otimes 1)\delta(x)+(1\otimes \ad_x)\Delta(y)+(1\otimes \ad_y)\tau\big(\Delta(x)\big)&=&\delta(x\cdot_P y).
\end{eqnarray}
\end{enumerate}
A noncommutative pseudo-Poisson bialgebra $((P,\cdot_P,\Poisson{-,-}_P),\Delta,\delta)$ is called a {\bf noncommutative Poisson bialgebra} if $P^*$ is also a coherent noncommutative Poisson algebra.
 \end{defi}

 \begin{rmk}
Even though $(P;L,R,\ad)$ is a representation of a noncommutative Poisson algebra $(P,\cdot_P,\Poisson{-,-}_P) $, $(P\otimes P;1\otimes L,R\otimes 1,\ad\otimes 1+1\otimes \ad)$ is just a quasi-representation. This  is the reason why we need to use quasi-representations to deal with the bialgebra theory of noncommutative Poisson algebras.
 \end{rmk}
\begin{rmk}
   Let $((P,\cdot_P,\Poisson{-,-}_P),\Delta,\delta)$ be a noncommutative Poisson bialgebra. Then by the fact that $(P,\delta)$ is a Lie coalgebra and \eqref{eq:PBA 3}, $((P,\Poisson{-,-}_P),\delta)$ is a Lie bialgebra. By the fact that $(P,\Delta)$ is a coassociative coalgebra, \eqref{eq:PBA 1} and \eqref{eq:PBA 4}, $((P,\cdot_P),\Delta)$ is an  antisymmetric infinitesimal bialgebra.
\end{rmk}

\emptycomment{Furthermore, let $\theta_\ast:P^*\longrightarrow P^*\wedge P^*$ and $\Theta_\ast: P^*\longrightarrow P^*\otimes P^*$ be the linear maps induced by the operations $\{-,-\}_P$ and $\cdot_P$, respectively. Then \eqref{eq:PBA 1}, \eqref{eq:PBA 3} and \eqref{eq:PBA 5} are equivalent to that $\theta_\ast+\Theta_\ast$ is a $1$-cocycle on the noncommutative Poisson algebra $(P^*,\delta^*,\Delta^*)$ with   coefficients in the quasi-representation $(P^*;1\otimes \frkL,\frkR\otimes 1,\frad)$. \yh{ I think there is mistakes. should it be ``$\theta_\ast+\Theta_\ast$ is a $1$-cocycle on the noncommutative Poisson algebra $(P^*,\delta^*,\Delta^*)$ with   coefficients  in the  representation $(P^*\otimes P^*;1\otimes \frkL,\frkR\otimes 1,\frad)$???''}}

\begin{thm}\label{thm:equivalent of Poisson bialgebra}
  Let $(P,\cdot_P,\Poisson{-,-}_P)$ be a coherent noncommutative Poisson algebra with two comultiplications $\Delta:P\longrightarrow \otimes^2P$ and $\delta:P\longrightarrow \wedge^2P$. Suppose that $\delta^*$ and $\Delta^*$ induce a coherent noncommutative Poisson algebra structure on $P^*$. Set $\Poisson{\alpha,\beta}_{P^*}=\delta^*(\alpha,\beta)$ and $\alpha\cdot_{P^*}\beta=\Delta^*(\alpha,\beta)$. Then the following statements are  equivalent:
   \begin{enumerate}
\item[\rm(i)]$((P,\cdot_P,\Poisson{-,-}_P),\Delta,\delta)$ is a noncommutative Poisson bialgebra.
\item[\rm(ii)]$(P,P^*;-R^*,-L^*,\ad^*,-\frkR^*,-\frkL^*,\frad^*)$ is a matched pair of coherent noncommutative Poisson algebras, where $\frkR^*,\frkL^*$ and $\frad^*$ are given by
 \begin{equation}\label{eq:frkRLad}
  \langle \frkR^*_\alpha x,\beta\rangle=-\langle x,\beta\cdot_{P^*}\alpha\rangle,\quad  \langle \frkL^*_\alpha x,\beta\rangle=-\langle x,\alpha\cdot_{P^*}\beta\rangle,\quad  \langle \frad^*_\alpha x,\beta\rangle=-\langle x,\{\alpha,\beta\}_{P^*}\rangle
 \end{equation}
 for all $x\in P,\alpha,\beta\in P^*$. \emptycomment{Furthermore, the noncommutative Poisson algebra structure on $P\oplus P^*$ is coherent.\yh{The last sentence holds naturally? Or this is another condition? Namely, is it enough to require that $(P,P^*;-R^*,-L^*,\ad^*,-\frkR^*,-\frkL^*,\frad^*)$ is a matched pair of coherent noncommutative Poisson algebras?}\jf{The last sentence holds naturally. The meaning of this sentence has been contained in the following statement (iii). I think this sentence should be deleted.}}
\item[\rm(iii)]$((P\oplus P^*,\frkB),P,P^*)$ is a Manin triple of noncommutative Poisson algebras, where the invariant symmetric bilinear form $\frkB$ on $P\oplus P^*$ is given by
\begin{equation}\label{eq:MTP-natural bilinear form}
  \frkB(x+\alpha,y+\beta)=\langle x,\beta\rangle+\langle \alpha,y\rangle,\quad\forall~x,y\in P,\alpha,\beta\in P^*.
\end{equation}
\end{enumerate}
\end{thm}
\pf First, we show that (i) and (ii) are equivalent. It is known that $((P,\Poisson{-,-}_P),\delta_P)$ is a Lie bialgebra if and only if $(P,P^*;\ad^*,\frad^*)$ is a matched pair of Lie algebras and $((P,\cdot_P),\Delta_P)$ is an antisymmetric infinitesimal bialgebra if and only if $(P,P^*;-R^*,-L^*,-\frkR^*,-\frkL^*)$ is a matched pair of associative algebras. By a straightforward calculation, we have
$$\eqref{eq:PBA 2}\Longleftrightarrow\eqref{eq:MPP1}\Longleftrightarrow\eqref{eq:MPP4} \quad \mbox{and}
\quad \eqref{eq:PBA 5}\Longleftrightarrow\eqref{eq:MPP2}\Longleftrightarrow\eqref{eq:MPP3}, $$
in which $\huaL=L^*,\huaR=R^*,\rho=\ad^*$ and $\tilde{\huaL}=\frkL^*,\tilde{\huaR}=\frkR^*,\varrho=\frad^*$.

Next, we show that (ii) and (iii) are equivalent. It is known  (\cite{CP})  that $(P,P^*;\ad^*,\frad^*)$ is a matched pair of Lie algebras if and only if $(P\oplus P^*,P,P^*;\frkB)$ is a Manin triple of Lie algebras. Similarly, $(P, P^*; -R^*, -L^*, -\frkR^*, -\frkL^*)$ is a matched pair of associative algebras if and only if $(P\oplus P^*,P,P^*;\frkB)$  is a double construction of Frobenius algebra (\cite{Bai1}, also see Remark~\ref{rmk:Manin}). The rest conditions in a matched pair of coherent noncommutative Poisson algebras are equivalent to the Leibniz rule that the noncommutative Poisson algebra structure on $P\oplus P^*$ satisfies. Thus, (ii) and (iii) are equivalent.\qed

\begin{rmk}
 Recall (\cite{WB}) that a {\bf compatible Lie bialgebra} structure
on a compatible Lie algebra $(\mathfrak{g},[-,-]_{1},[-,-]_{2})$ is a
pair of linear maps $\alpha, \beta:\mathfrak{g}\rightarrow
\mathfrak{g}\otimes\mathfrak{g}$ such that for any $k_1,k_2\in
\mathbb{K},$ $k_1\alpha+k_2\beta$ is a Lie bialgebra structure on
the Lie algebra $(\mathfrak{g},k_1[-,-]_{1}+k_2[-,-]_{2})$. A compatible Lie bialgebra is equivalent to a Manin triple of
compatible Lie algebras in the following sense: assume that $(\mathfrak{g},[-,-]_{1},[-,-]_{2})$ and $(\mathfrak{g}^*,\{-,-\}_{1},\{-,-\}_{2})$
are compatible Lie algebras, there is a compatible Lie algebra structure $(\Courant{-,-}_1, \Courant{-,-}_2)$ on the direct sum of the underlying vector spaces of $\mathfrak{g}$ and  $\mathfrak{g}^*$ such that
the bilinear form given by \eqref{eq:MTP-natural bilinear form} is invariant in the sense
$$\frkB(\Courant{u,v}_1,w)=\frkB(u, \Courant{v,w}_1), \frkB(\Courant{u,v}_2,w)=\frkB(u, \Courant{v,w}_2),\quad\forall~u,v,w\in \g\oplus\g^*.$$
It is obvious that a Manin triple of noncommutative Poisson algebras naturally gives a Manin triple of compatible Lie algebras and hence a noncommutative Poisson bialgebra naturally gives
a compatible Lie bialgebra.

\end{rmk}

\section{Coboundary (pseudo-)Poisson bialgebras}\label{sec:coboundary}
To begin with, we recall some important results of   coboundary Lie bialgebras and   coboundary antisymmetric
infinitesimal bialgebras.

Let $(\g,[-,-]_\g)$ be a Lie algebra and $r\in \g\otimes\g$. Then the linear map $\delta$ defined by \eqref{eq:coboundary Lie cobracket} makes $(\g,\delta)$ into a coboundary  Lie bialgebra if and only if the following conditions are satisfied
\begin{enumerate}
\item[\rm(i)]$(\ad_x\otimes 1 +1\otimes \ad_x)(r+\tau(r))=0$,
\item[\rm(ii)]$(\ad_x\otimes 1\otimes 1 +1\otimes\ad_x\otimes 1+1\otimes 1\otimes \ad_x)([r_{12},r_{13}]+[r_{12},r_{23}]+[r_{13},r_{23}])=0$.
\end{enumerate}
In particular, the following equation
\begin{equation}
  {\bf C}(r)=[r_{12},r_{13}]+[r_{12},r_{23}]+[r_{13},r_{23}]=0
\end{equation}
is called the {\bf classical Yang-Baxter equation ($\CYBE$)}.

Let $(A,\cdot)$ be an associative algebra and $r\in \g\otimes\g$.
Then the linear map $\Delta$ defined by \eqref{eq:coboundary
associative cobracket} makes $(A,\Delta)$ into a coboundary
 antisymmetric
infinitesimal bialgebra if and only if the following conditions
are satisfied
\begin{enumerate}
\item[\rm(i)]$(L_x\otimes 1-1\otimes R_x)(1\otimes L_y-R_y\otimes 1)(r+\tau(r))=0$,
\item[\rm(ii)]$(1\otimes 1\otimes L_x-R_x\otimes 1\otimes1)(r_{12}\cdot r_{13}+r_{13}\cdot r_{23}-r_{23}\cdot r_{12})=0$.
\end{enumerate}
In particular, the following equation
\begin{equation}
  {\bf A}(r)=r_{12}\cdot r_{13}+r_{13}\cdot r_{23}-r_{23}\cdot r_{12}=0
\end{equation}
is called the {\bf associative Yang-Baxter equation ($\AYBE$)}. See \cite{A1,A3,Bai1,BGN2013A} for more details.

Next, we introduce the notion of a coboundary noncommutative (pseudo-)Poisson bialgebra.
\begin{defi}
  A noncommutative (pseudo-)Poisson bialgebra $(P,\Delta,\delta)$ is called {\bf coboundary} if there exists an $r\in P\otimes P$ such that
  \begin{equation}\label{eq:coboundary Poisson}
\delta(x)=(1\otimes \ad_x+\ad_x\otimes 1) r,\quad \Delta(x)=(1\otimes L_x -R_x\otimes 1 )r,\quad\forall~x\in P.
\end{equation}
\end{defi}
We denote a coboundary noncommutative (pseudo-)Poisson bialgebra  by $(P,r)$.
\begin{thm}\label{thm:coboundary pseudo-Poisson bialgebra}
Let $(P,\cdot_P,\Poisson{-,-}_P)$ be a coherent noncommutative Poisson algebra and $r\in P\otimes P$. Define $\Delta:P\longrightarrow P\otimes P$ and $\delta:P\longrightarrow \wedge^2P$ by \eqref{eq:coboundary Poisson}. Then $(P,\Delta,\delta)$ is a noncommutative pseudo-Poisson bialgebra if and only if the following conditions are satisfied:
\begin{enumerate}
\item[\rm(i)]$(L_x\otimes 1-1\otimes R_x)(1\otimes L_y-R_y\otimes 1)(r+\tau(r))=0$;
\item[\rm(ii)]$(\ad_x\otimes 1 +1\otimes \ad_x)(r+\tau(r))=0;$
\item[\rm(iii)]$( \ad_x\otimes 1)(1\otimes L_y-R_y\otimes 1)(r+\tau(r))=0$;
\item[\rm(iv)]$(1\otimes 1\otimes L_x-R_x\otimes 1\otimes1){\bf A}(r)=0$;
\item[\rm(v)]$(\ad_x\otimes 1\otimes 1 +1\otimes\ad_x\otimes 1+1\otimes 1\otimes \ad_x){\bf C}(r)=0$;
\item[\rm(vi)]for $r=\sum_ja_j\otimes b_j\in P\otimes P$,
\begin{eqnarray*}
&& (\ad_x\otimes 1\otimes 1){\bf A}(r)+(1\otimes1\otimes L_x-1\otimes R_x\otimes 1){\bf C}(r)\\
&&+\sum_j\tau\Big((1\otimes\ad_{a_j})(1\otimes L_x-R_x\otimes 1)(r+\tau(r))\Big)\otimes b_j=0.
\end{eqnarray*}
\end{enumerate}
\end{thm}
\pf First by Conditions (i) and (iv),  the dual map $ \Delta^*:\otimes^2 P^*\longrightarrow P^*$ defines an associative algebra structure on $P^*$. By Conditions (ii) and (v),  the dual map
$\delta^*:\wedge^2 P^*\longrightarrow P^*$ defines a Lie algebra structure on $P^*$. By the fact that the noncommutative Poisson algebra $P$ is coherent and (vi), we can deduce that the Leibniz rule is satisfied. Thus, $(P^*,\Delta^*,\delta^*)$ is a noncommutative Poisson algebra, which implies that Condition (a) in Definition \ref{defi:quasiPBA} holds.

Then since $\Delta+\delta$ is a coboundary, it is naturally closed. This implies that Condition (b) in Definition \ref{defi:quasiPBA} holds.

Finally, by the definition of $\delta$ given by \eqref{eq:coboundary Poisson}, it is straightforward to deduce that \eqref{eq:PBA 4} holds. By a direct calculation, we can obtain that \eqref{eq:PBA 5} is equivalent to
$$\big(1\otimes(\ad_{x\cdot_P y}-\ad_x L_y)\big)r+(1\otimes \ad_y R_x)\tau(r)-\big(L_x\otimes \ad_y\big)(r+\tau(r))=0.$$
Since the noncommutative Poisson algebra $(P,\cdot_P,\Poisson{-,-}_P)$ is coherent, we have
$$\ad_{x\cdot_P y}-\ad_x L_y=\ad_y R_x.$$
Thus \eqref{eq:PBA 5} is equivalent to
$$(1\otimes\ad_y R_x)(r+\tau(r))-\big(L_x\otimes \ad_y\big)(r+\tau(r))=0,$$
which is equivalent to Condition (iii). Thus Condition (c) in Definition \ref{defi:quasiPBA} holds.

The converse can be proved similarly. We omit the details.\qed\vspace{3mm}

\begin{defi}
 Let $(P,\cdot_P,\Poisson{-,-}_P)$ be a coherent
 noncommutative  Poisson algebra and $r\in P\otimes P$. The equation
 \begin{equation}
   {\bf A}(r)={\bf C}(r)=0
 \end{equation}
 is called the {\bf Poisson Yang-Baxter equation ($\PYBE$)} in $P$.
\end{defi}

Let $(\g,[-,-]_\g)$ be a Lie algebra and $r\in\wedge^2\g$ a solution of the $\CYBE$. Then the Lie algebra structure $[-,-]_{\g^*}:\g^*\times\g^*\longrightarrow\g^*$ induced by $r$ is given by
\begin{equation}\label{eq:Lie coalgebra}
  [\alpha,\beta]_{\g^*}=\ad^*_{r^\sharp(\alpha)}\beta-\ad^*_{r^\sharp(\beta)}\alpha,\quad~\forall~\alpha,\beta\in \g^*,
\end{equation}
where $r^\sharp:\g^*\longrightarrow \g$ is defined by
\begin{equation}
  \langle r^\sharp(\alpha),\beta\rangle=r(\alpha,\beta).
\end{equation}
Furthermore, $r^\sharp$ is a Lie algebra homomorphism from the Lie algebra $(\g^*,[-,-]_{\g^*})$ to $(\g,[-,-]_\g)$.

Let $(A,\cdot_A)$ be an associative algebra and $r\in\wedge^2 A$ a solution of the $\AYBE$. Then the associative algebra structure $\cdot_{A^*}:A^*\otimes A^*\longrightarrow A^*$ induced by $r$ is given by
\begin{equation}\label{eq:associative coalgebra}
 \alpha\cdot_{A^*}\beta=-R^*_{r^\sharp(\alpha)}\beta-L^*_{r^\sharp(\beta)}\alpha,\quad~\forall~\alpha,\beta\in A^*.
\end{equation}
Furthermore, $r^\sharp$ is an associative algebra homomorphism from the associative algebra $(A^*,\cdot_{A^*})$ to $(A,\cdot_A)$.

The following theorem shows that a skew-symmetric solution of the $\PYBE$ in a coherent noncommutative Poisson algebra gives rise to a noncommutative Poisson bialgebra.

\begin{thm}\label{thm:Poisson r-matrix}
 Let $(P,\cdot_P,\Poisson{-,-}_P)$ be a coherent noncommutative Poisson algebra and $r\in \wedge^2 P$ a solution of the $\PYBE$. Then the maps $\cdot_{P^*}:=\Delta^*:\otimes^2 P^*\longrightarrow P^*$ and $\Poisson{-,-}_{P^*}:=\delta^*:\wedge^2 P^*\longrightarrow P^*$, where $\Delta$ and $\delta$ are given by \eqref{eq:coboundary Poisson}, induce a coherent noncommutative Poisson algebra structure on $P^*$ such that $(P,P^*)$ is a noncommutative Poisson bialgebra.
\end{thm}
\pf By the fact that $r$ is skew-symmetric, we deduce that (i)-(iii) in Theorem \ref{thm:coboundary pseudo-Poisson bialgebra} hold.  By the fact that $r$ is a skew-symmetric solution of the $\PYBE$,  we deduce that (iv)-(vi) in Theorem \ref{thm:coboundary pseudo-Poisson bialgebra} hold.   Thus, $(P,P^*)$ is a noncommutative pseudo-Poisson bialgebra.

 Then  by   \eqref{eq:Lie coalgebra} and \eqref{eq:associative coalgebra}, we have
\begin{eqnarray*}
  &&\Poisson{\alpha\cdot_{P^*} \beta,\gamma}_{P^*}-\Poisson{\alpha,\beta\cdot_{P^*} \gamma}_{P^*}-\Poisson{\beta,\gamma\cdot_{P^*} \alpha}_{P^*}\\
  &=&\ad^*_{r^\sharp(\alpha\cdot_{P^*} \beta)}\gamma-\ad^*_{r^\sharp(\gamma)}(\alpha\cdot_{P^*} \beta)+\ad^*_{r^\sharp(\beta\cdot_{P^*} \gamma)}\alpha-\ad^*_{r^\sharp(\alpha)}(\beta\cdot_{P^*} \gamma)\\
  &&+\ad^*_{r^\sharp(\gamma\cdot_{P^*} \alpha)}\beta-\ad^*_{r^\sharp(\beta)}(\gamma\cdot_{P^*} \alpha)\\
  &=&\ad^*_{r^\sharp(\alpha)\cdot_{P} r^\sharp(\beta)}\gamma+\ad^*_{r^\sharp(\gamma)}R^*_{r^\sharp(\alpha)}\beta+\ad^*_{r^\sharp(\gamma)}L^*_{r^\sharp(\beta)}\alpha
  +\ad^*_{r^\sharp(\beta)\cdot_{P} r^\sharp(\gamma)}\alpha\\
  &&+\ad^*_{r^\sharp(\alpha)}R^*_{r^\sharp(\beta)}\gamma+\ad^*_{r^\sharp(\alpha)}L^*_{r^\sharp(\gamma)}\beta+\ad^*_{r^\sharp(\gamma)\cdot_{P} r^\sharp(\alpha)}\beta+\ad^*_{r^\sharp(\beta)}R^*_{r^\sharp(\gamma)}\alpha+\ad^*_{r^\sharp(\beta)}L^*_{r^\sharp(\alpha)}\gamma\\
  &=&\big(\ad^*_{r^\sharp(\alpha)\cdot_{P} r^\sharp(\beta)}+\ad^*_{r^\sharp(\alpha)}R^*_{r^\sharp(\beta)}+\ad^*_{r^\sharp(\beta)}L^*_{r^\sharp(\alpha)}\big)\gamma+\big(\ad^*_{r^\sharp(\gamma)\cdot_{P} r^\sharp(\alpha)}+\ad^*_{r^\sharp(\gamma)}R^*_{r^\sharp(\alpha)}\\
  &&+\ad^*_{r^\sharp(\alpha)}L^*_{r^\sharp(\gamma)}\big)\beta+\big(\ad^*_{r^\sharp(\beta)\cdot_{P} r^\sharp(\gamma)}+\ad^*_{r^\sharp(\beta)}R^*_{r^\sharp(\gamma)}+\ad^*_{r^\sharp(\gamma)}L^*_{r^\sharp(\beta)}\big)\alpha\\
  &=&0.
\end{eqnarray*}
The last equality holds because
$$\langle\ad^*_{x\cdot_P y}\alpha+\ad^*_x R^*_y\alpha +\ad^*_y L^*_x\alpha,z\rangle=\langle x\cdot_P\{y,z\}_P+\{x,z\}_P\cdot_P y-\{x\cdot_P y, z\}_P,\alpha \rangle=0,$$
which follows from the Leibniz rule of the noncommutative Poisson algebra $P$. We deduce that $(P^*,\cdot_{P^*},\Poisson{-,-}_{P^*})$ is a coherent noncommutative Poisson algebra. Thus, $(P,P^*)$ is a noncommutative Poisson bialgebra. \qed

\begin{rmk}
 The study of coboundary compatible Lie bialgebras and the classical Yang-Baxter equation in compatible Lie algebras was also given in
\cite{WB}. It is straightforward to see that a skew-symmetric solution of the Poisson Yang-Baxter equation in a coherent noncommutative Poisson algebra $(P,\cdot_P,\{-,-\}_P)$
 is a skew-symmetric solution of
the classical Yang-Baxter equation in the compatible Lie algebra $(P, \{-,-\}_P,[-,-]_P)$ and hence gives a (coboundary) compatible Lie bialgebra.

\end{rmk}

\begin{cor}
  Let $(P,\cdot_P,\Poisson{-,-}_P)$ be a coherent noncommutative Poisson algebra and $r\in \wedge^2 P$ a solution of the $\PYBE$. Then $r^\sharp:P^*\longrightarrow P$ is a noncommutative Poisson algebra homomorphism from the coherent noncommutative Poisson algebra $(P^*,\cdot_{P^*},\Poisson{-,-}_{P^*})$ to $(P,\cdot_P,\Poisson{-,-}_P)$.
\end{cor}

\begin{ex}{\rm
 Let $(A,\cdot_A,\{-,-\}_\hbar)$ be the standard noncommutative Poisson algebra given by Example \ref{ex:coherent poisson 2}. If $r\in\wedge^2A$ is a solution of the $\AYBE$ in the associative algebra $(A,\cdot_A)$, then $r$ is a solution of the $\PYBE$  in the standard noncommutative Poisson algebra $A$. Furthermore, the dual noncommutative Poisson algebra $(A^*,\cdot_{A^*},\Poisson{-,-}_{P^*})$ induced by $r$ is also standard, where
\begin{eqnarray*}
 \alpha\cdot_{A^*}\beta&=&-R^*_{r^\sharp(\alpha)}\beta-L^*_{r^\sharp(\beta)}\alpha,\\
 \Poisson{\alpha,\beta}_{P^*}&=& \hbar(\alpha\cdot_{A^*}\beta- \beta\cdot_{A^*}\alpha),\quad\forall~\alpha,\beta\in A^*.
 \end{eqnarray*}
 }
\end{ex}

\begin{ex}{\rm
Consider the coherent noncommutative Poisson algebra $P$
given by Example \ref{ex:SNCPA1}, then $  r$   given
by
$$  r =\kappa_{13} e_1\wedge e_3+\kappa_{23} e_2\wedge e_3$$
is a solution of the $\PYBE$ in $P$, where $\kappa_{13}$ and $\kappa_{23}$ are constants.
}

\end{ex}
\begin{ex}{\rm
Consider the coherent noncommutative Poisson algebra $P$
given by Example \ref{ex:SNCPA2}, then $ r$ given
by
$$  r=\kappa_{12} e_1\wedge e_2+\kappa_{14} e_1\wedge e_4+\kappa_{12} e_2\wedge e_3+\kappa_{24}e_2\wedge e_4-\kappa_{14}e_3\wedge e_4$$
and
$$ r=\kappa_{14} e_1\wedge e_4+\kappa_{24} e_2\wedge e_4+(\kappa_{14}+\kappa_{24})e_3\wedge e_4$$
are solutions of the $\PYBE$ in $P$, where $\kappa_{12}$,
$\kappa_{14}$ and $\kappa_{24}$ are constants. }

\end{ex}
Let $(P,\cdot_P,\Poisson{-,-}_P)$ be a coherent
noncommutative Poisson algebra and $r\in \wedge^2 P$. Assume that
$r$ is nondegenerate, i.e.
$r^\sharp:P^*\longrightarrow P$ is invertible. Define $\omega\in
\wedge^2 P^*$ by
$$\omega(x,y)=\langle (r^\sharp)^{-1}(x),y\rangle,\quad\forall~x,y\in P.$$
Then we have
\begin{pro}\label{pro:cocycle and symplectic}
With the above notations, $r\in \wedge^2 P$ is a solution of the $\PYBE$  in a coherent noncommutative Poisson algebra $(P,\cdot_P,\Poisson{-,-}_P)$ if and only if $\omega$ is both a  Connes cocycle on the associative $(P,\cdot_P)$ and a symplectic structure on the Lie algebra $(P,\Poisson{-,-}_P)$, i.e.,
 \begin{eqnarray*}
   \omega(x\cdot_P y,z)+\omega(y\cdot_P z,x)+\omega(z\cdot_P x,y)&=&0,\\
    \omega(\Poisson{x, y}_P,z)+\omega(\Poisson{y, z}_P,x)+\omega(\Poisson{z,x}_P,y)&=&0,\quad\forall~x,y,z\in P.
 \end{eqnarray*}
\end{pro}
\pf It follows from the fact that $r$ is a solution of the $\AYBE$  if and only if $\omega$ is a Connes cocycle on the associative $(P,\cdot_P,)$ (\cite{Bai1})
 and
$r$ is a solution of the $\CYBE$  if and only if $\omega$ is a symplectic structure on the Lie algebra $(P,\Poisson{-,-}_P)$ (\cite{Drinfeld}).
\qed\vspace{3mm}

Let $(P,\cdot_P,\Poisson{-,-}_P)$ be a coherent noncommutative Poisson algebra. An element $\frks\in P\otimes P$ is called   {\bf $(L,R,\ad)$-invariant} if
\begin{equation}\label{eq:LRad-invariant}
 (1\otimes L_x -R_x\otimes 1)\frks=0,\quad (\ad_x\otimes 1 +1\otimes \ad_x)\frks=0,\quad\forall~x\in P.
 \end{equation}

\begin{pro}\label{pro:symmetric part PYBE}
   Let $(P,\cdot_P,\Poisson{-,-}_P)$ be a coherent noncommutative Poisson algebra. Let $r=\frka+\frks\in P\otimes P$ with skew-symmetric part $\frka$ and symmetric part $\frks$. If the symmetric part $\frks$ of $r$ is $(L,R,\ad)$-invariant and $r$ is a solution of the $\PYBE$, then $(P,r)$ gives a coboundary noncommutative pseudo-Poisson bialgebra. Furthermore, this coboundary noncommutative pseudo-Poisson bialgebra is a noncommutative Poisson bialgebra if and only if the symmetric part $\frks$ of $r$ satisfies
  \begin{equation}\label{eq:symmetric conditon}
    \ad^*_{\frks^\sharp(\alpha)\cdot_P\frks^\sharp(\beta) }\gamma +  \ad^*_{\frks^\sharp(\gamma)\cdot_P\frks^\sharp(\alpha) }\beta+  \ad^*_{\frks^\sharp(\beta)\cdot_P\frks^\sharp(\gamma) }\alpha=0,
  \end{equation}
  where $\frks^\sharp:P^*\longrightarrow P$ is defined by
  \begin{equation}
    \langle\frks^\sharp(\alpha),\beta\rangle=\frks(\alpha,\beta),\quad\forall~\alpha,\beta\in P^*.
  \end{equation}
\end{pro}
\pf Since the symmetric part $\frks$ of $r$ is $(L,R,\ad)$-invariant and $r$ is a solution of the $\PYBE$, by Theorem \ref{thm:coboundary pseudo-Poisson bialgebra}, $(P,r)$ gives a coboundary noncommutative pseudo-Poisson bialgebra.

By the invariance of $\frks$ and ${\bf A}(r)=0$, we have
$${\frka^\sharp(\alpha)\cdot_{P} \frka^\sharp(\beta)}-\frka^\sharp(\alpha\cdot_{P^*}\beta)=\frks^\sharp(\alpha)\cdot_P\frks^\sharp(\beta),$$
where $\alpha\cdot_{P^*}\beta=-R^*_{\frka^\sharp(\alpha)}\beta-L^*_{\frka^\sharp(\beta)}\alpha.$ By this fact, with a similar proof of Theorem \ref{thm:Poisson r-matrix}, the noncommutative Poisson algebra $(P^*,\cdot_{P^*},\{-,-\}_{P^*})$ is coherent with $ \{\alpha,\beta\}_{P^*}=\ad^*_{\frka^\sharp(\alpha)}\beta-\ad^*_{\frka^\sharp(\beta)}\alpha$ if and only if the symmetric part $\frks$ of $r$ satisfies \eqref{eq:symmetric conditon}. We omit the details.\qed\vspace{3mm}

At the end of this section, we establish the Drinfeld  double theory for  noncommutative Poisson bialgebras.
Let $((P,\cdot_P,\Poisson{-,-}_P),\Delta,\delta)$ be a noncommutative Poisson bialgebra. By Theorem \ref{thm:Manin-Standard} and Theorem \ref{thm:equivalent of Poisson bialgebra}, there is a coherent noncommutative Poisson algebra structure on $\huaD=P\oplus P^*$ with the associative multiplication $\ast_\huaD:\huaD\times \huaD\longrightarrow \huaD$ and the Lie bracket $\{-,-\}_{\huaD}:\huaD\times \huaD\longrightarrow \huaD$ given by
 \begin{eqnarray}
  (x+\alpha)\ast_\huaD (y+\beta)&=&x\cdot_P y-\frkR^*_\alpha y-\frkL^*_\beta x+\alpha\cdot_{P^*}\beta-R^*_x \beta-L^*_y\alpha,\\
  \Poisson{x+\alpha,y+\beta}_{\huaD}&=&\Poisson{x,y}_{P}+\frad^*_{\alpha}y-\frad^*_{\beta}x+\Poisson{\alpha,\beta}_{P^*}+\ad^*_{x}\beta-\ad^*_{y}\alpha
  \end{eqnarray}
for all $x,y\in P,\alpha,\beta\in P^*$, where $\frkR^*,\frkL^*$ and $\frad^*$ are given by \eqref{eq:frkRLad}.

Let $\{e_1,\cdots,e_n\}$ be a basis of $P$ and $\{e^*_1,\cdots,e^*_n\}$ its dual basis. Let $r=\sum_ie_i\otimes e^*_i\in\huaD\otimes \huaD$.

\begin{pro}
 With the above notations,  $r=\sum_ie_i\otimes e^*_i\in\huaD\otimes \huaD$ is a solution of the $\PYBE$ in the coherent noncommutative Poisson algebra $\huaD$ such that $(\huaD,\huaD^*)$ is a noncommutative Poisson bialgebra.
\end{pro}
\pf Obviously the symmetric part of $r$ is $\frks=\half\sum_i(e_i\otimes e^*_i+e^*_i\otimes e_i)$ and the skew-symmetric part of $r$ is $\frka=\half\sum_i(e_i\otimes e^*_i-e^*_i\otimes e_i)$. By the Drinfeld double theory of Lie algebras and the associative double theory of associative algebras (see \cite{Bai1} for more details on the associative double theory), we know that $\frks$ is $(L,R,\ad)$-invariant and $r$ satisfies the classical and the associative Yang-Baxter equations. Thus $r$ is a solution of the $\PYBE$, and by Proposition \ref{pro:symmetric part PYBE}, $(\huaD,r)$ is a noncommutative pseudo-Poisson bialgebra. By a straightforward  calculation,   the noncommutative Poisson algebra structure on $\huaD^*$ is given by
\begin{eqnarray*}
  (x+\alpha)\ast_{\huaD^*} (y+\beta)&=&x\cdot_P y-\alpha\cdot_{P^*}\beta,\\
    \Poisson{x+\alpha,y+\beta}_{\huaD^*}&=&\Poisson{x,y}_{P}-\Poisson{\alpha,\beta}_{P^*}.
  \end{eqnarray*}
It is easy to check that the noncommutative Poisson algebra structure on $\huaD^*$ is coherent. Thus $(\huaD,\huaD^*)$ is a noncommutative Poisson bialgebra.\qed

\section{Noncommutative  pre-Poisson algebras, Rota-Baxter operators and  $\huaO$-operators}\label{sec:pre-Poisson}

In this section, we introduce the notions of  noncommutative  pre-Poisson algebras and Rota-Baxter operators (more generally $\huaO$-operators) on a noncommutative Poisson algebra. We show that on the one hand,
an $\huaO$-operator on a noncommutative Poisson algebra gives
a noncommutative  pre-Poisson algebra, and on the other hand, a noncommutative  pre-Poisson algebra naturally gives an $\huaO$-operator on the sub-adjacent noncommutative Poisson algebra.
 We use $\huaO$-operators and noncommutative  pre-Poisson algebras to construct some skew-symmetric solutions of the Poisson Yang-Baxter equation.

\subsection{Noncommutative  pre-Poisson algebras}

Recall that a {\bf pre-Lie algebra} is a pair $(A,\ast)$, where
$A$ is a vector space, and  $\ast:A\otimes A\longrightarrow A$ is
a bilinear multiplication satisfying that for all $x,y,z\in A$,
the associator $(x,y,z)=(x\ast y)\ast z-x\ast(y\ast z)$ is
symmetric in $x,y$, i.e.,
$$(x,y,z)=(y,x,z),\;\;{\rm or}\;\;{\rm
equivalently,}\;\;(x\ast y)\ast z-x\ast(y\ast z)=(y\ast x)\ast
z-y\ast(x\ast z).$$

\begin{lem}{\rm(\cite{Pre-lie})}\label{lem:pre-Lie-Lie} Let $(A,\ast)$ be a pre-Lie algebra. The commutator
$ [x,y]_A=x\ast y-y\ast x$ defines a Lie algebra structure on $A$,
which is called the {\bf sub-adjacent Lie algebra} of $(A,\ast)$ and denoted by $A^c$. Furthermore, $L:A\rightarrow
\gl(A)$ defined by
\begin{equation}\label{eq:defiLpreLie}
L_xy=x\ast y,\quad \forall x,y\in A
\end{equation}
 gives a representation of $A^c$ on $A$.
\end{lem}

There is a similar relationship between dendriform algebras and associative algebras.

\begin{lem}{\rm(\cite{Lod2})}\label{lem:den-ass}
Let $(A,\succ,\prec)$ be a dendriform algebra. Then $(A,\cdot)$ is an associative algebra, where $x\cdot y=x\succ y+x\prec y$. Moreover, for $x\in A$, define $L_{\succ x},R_{\prec x}:A\longrightarrow\gl(A)$ by
\begin{equation}\label{eq:dendriform-rep}
L_{\succ x}(y)=x\succ y,\quad R_{\prec x}(y)=y\prec x,\quad\forall~y\in A.
\end{equation}
Then $(A;L_\succ,R_\prec)$ is a representation of the associative algebra $(A,\cdot)$.
\end{lem}

Now we are ready to give the main notion in this subsection.

\begin{defi}
  A {\bf noncommutative  pre-Poisson algebra} is a quadruple $(A,\succ,\prec,\ast)$ such that $(A,\succ,\prec)$ is a dendriform algebra and $(A,\ast)$ is a pre-Lie algebra satisfying the following compatibility conditions:
  \begin{eqnarray}
    \label{eq:pre-Poisson 1}(x\ast y-y\ast x)\succ z&=&x\ast(y\succ z)-y\succ(x\ast z),\\
    \label{eq:pre-Poisson 2} x\prec(y\ast z-z\ast y)&=&y\ast(x\prec z)-(y\ast x)\prec z,\\
    \label{eq:pre-Poisson 3} (x\succ y+x\prec y)\ast z&=&(x\ast z)\prec y+x\succ(y\ast z).
  \end{eqnarray}
  A noncommutative  pre-Poisson algebra $(A,\succ,\prec,\ast)$  is called {\bf coherent} if it also satisfies
  \begin{equation}\label{eq:pre-Poisson coherent}
    (x\succ y+x\prec y)\ast z=x\ast(y\succ z)+y\ast(z\prec x).
  \end{equation}
\end{defi}

\begin{rmk}
Aguiar introduced the  notion of pre-Poisson algebras in  \cite{A2}, as the splitting of Poisson algebras. If the dendriform algebra in a noncommutative pre-Poisson algebra reduces to a Zinbiel algebra, then we obtain a pre-Poisson algebra. Noncommutative pre-Poisson algebras can be viewed as the splitting of noncommutative Poisson algebras. See \cite{BBGN,PBG} for more details of splitting of operads.
\end{rmk}

\begin{pro}
 A noncommutative pre-Poisson algebra $(A,\succ,\prec,\ast)$ is  coherent if and only if it satisfies
 \begin{eqnarray*}
   (x\ast y)\circ z-x\circ(y\ast z)=(y\ast x)\circ z-y\circ(x\ast z)
 \end{eqnarray*}
 for all $x,y,z\in A$, where $x\circ y=x\succ y-y\prec x$.
\end{pro}
\pf By \eqref{eq:pre-Poisson 1}-\eqref{eq:pre-Poisson 3} in the definition of noncommutative pre-Poisson algebra , we have
\begin{eqnarray*}
 && x\ast(y\succ z)+y\ast(z\prec x)-(x\succ y+x\prec y)\ast z\\
   &=&(x\ast y-y\ast x)\succ z+y\succ(x\ast z)+z\prec(y\ast x-x\ast y)+(y\ast z)\prec x-(x\ast z)\prec y-x\succ(y\ast z)\\
   &=&\big((x\ast y)\succ z-z\prec(x\ast y)\big)-\big(x\succ(y\ast z)-(y\ast z)\prec x\big)-\big((y\ast x)\succ z-z\prec(y\ast x)\big)\\
   &&+\big(y\succ(x\ast z)-(x\ast z)\prec y\big)\\
   &=& (x\ast y)\circ z-x\circ(y\ast z)-(y\ast x)\circ z+y\circ(x\ast z).
\end{eqnarray*}
Then the conclusion follows immediately.\qed\vspace{3mm}

 Similarly, a coherent noncommutative pre-Poisson
algebra $(A,\succ,\prec,\ast)$ also satisfies
 $$(x\circ y)\ast z-x\ast(y\circ z)=(y\circ x)\ast z-y\ast(x\circ z).$$

\begin{defi}{ A {\bf compatible pre-Lie algebra} $(A,\circ,\ast)$
consists of two pre-Lie algebras $(A,\circ)$  and $(A,\ast)$ such that for any $k_1,k_2\in\mathbb{K}$, the following bilinear operation
\begin{equation}x\star y=k_1x\circ y+k_2x\ast y,\quad \forall~x,y\in A,\end{equation}
defines a pre-Lie algebra structure on $A$.
}\end{defi}

 It is straightforward to obtain the following result.

\begin{pro}
Let $(A,\circ)$  and $(A,\ast)$ be  two pre-Lie algebras. Then $(A,\circ,\ast)$ is a compatible pre-Lie algebra if and only if for all $ x,y,z\in A$,
\begin{equation}
(x\circ y)\ast z-x\ast (y\circ z)+(x\ast y)\circ z-x\circ (y\ast z)=(y\circ x)\ast z-y\ast (x\circ z)+(y\ast x)\circ z-y\circ (x\ast z).\end{equation}
\end{pro}

\begin{cor}
Let $(A,\succ,\prec,\ast)$ be a coherent noncommutative pre-Poisson algebra. Then $(A, \circ,\ast)$ is a compatible pre-Lie algebra, where $x\circ y=x\succ y-y\prec x$ for all $x,y\in A$.
\end{cor}

\begin{ex}\label{ex:coherent pre-Poisson}{\rm
  Let $(A,\succ,\prec)$ be a dendriform algebra. Then $(A,\succ,\prec,\ast_\hbar)$ is a coherent noncommutative  pre-Poisson algebra, where $\hbar$ is a fixed number and the pre-Lie algebra structure $\ast_\hbar$ is given by
   $$x\ast_\hbar y=\hbar(x\succ y-y\prec x),\quad \forall~x,y\in A.$$
  }
\end{ex}

A noncommutative pre-Poisson algebra gives rise to a noncommutative Poisson algebra and a representation on itself  naturally.

\begin{thm}
  Let $(A,\succ,\prec,\ast)$ be a noncommutative  pre-Poisson algebra.
  \begin{itemize}
  \item[{\rm(i)}]Define
  $$x\cdot y=x\succ y+x\prec y\quad\mbox{and}\quad \{x,y\}=x\ast y-y\ast x,\quad \forall~x,y\in A.$$
  Then $(A,\cdot,\{-,-\})$ is a noncommutative Poisson algebra, which is called the {\bf sub-adjacent noncommutative Poisson algebra} of $(A,\succ,\prec,\ast)$ and denoted by $A^c$.

  \item[{\rm(ii)}] If the noncommutative  pre-Poisson algebra $(A,\succ,\prec,\ast)$ is coherent, then $A^c$ is also coherent.

  \item[{\rm(iii)}] $(A;L_{\succ},R_\prec,L)$ is a representation   of the  sub-adjacent noncommutative Poisson algebra $A^c$, where $L_\succ$ and $R_\prec$ are given by \eqref{eq:dendriform-rep} and $L$ is given by \eqref{eq:defiLpreLie}.
      \end{itemize}
\end{thm}
\pf (i) By Lemma \ref{lem:pre-Lie-Lie} and \ref{lem:den-ass}, we deduce that $(A,\cdot)$ is an associative algebra and $(A,\{-,-\})$ is a Lie algebra.  By \eqref{eq:pre-Poisson 1}-\eqref{eq:pre-Poisson 3}, we have
\begin{eqnarray*}
\{x,y\cdot z\}-\{x,y\}\cdot z-y\cdot\{x,z\}&=&\Big(x\ast(y\succ z)-y\succ(x\ast z)-(x\ast y-y\ast x)\succ z\Big)\\
   &&+ \Big(x\ast(y\prec z)-(x\ast y)\prec z- y\prec(x\ast z-z\ast x)\Big)\\
  &&+\Big((y\ast x)\prec z+y\succ(z\ast x)-(y\succ z+y\prec z)\ast x\Big)\\
  &=&0,
\end{eqnarray*}
which implies that $(A,\cdot,\{-,-\})$ is a noncommutative Poisson algebra.

(ii) If the noncommutative pre-Poisson algebra $(A,\succ,\prec,\ast)$ is coherent, by \eqref{eq:pre-Poisson coherent}, we have
\begin{eqnarray*}
  \Poisson{x\cdot y,z}-\Poisson{x,y\cdot z}-\Poisson{y,z\cdot x}
  &=&(x\cdot y) \ast z-z\ast(x\succ y)-z\ast(x\prec y)+(y\cdot z) \ast x\\
  &&-x\ast(y\succ z)-x\ast(y\prec z)+(z\cdot x) \ast y-y\ast(z\succ x)-y\ast(z\prec x)\\
  &=&\big((x\cdot y) \ast z-x\ast(y\succ z)-y\ast(z\prec x)\big)+\big((z\cdot x) \ast y-z\ast(x\succ y)\\
  &&-x\ast(y\prec z)\big)+\big((y\cdot z) \ast x-y\ast(z\succ x)-z\ast(x\prec y)\big)\\
  &=&0,
\end{eqnarray*}
which implies that \eqref{eq:coboundary extra condition} holds. Thus, $(A,\cdot,\{-,-\})$ is coherent.

(iii) By Lemma \ref{lem:pre-Lie-Lie}, $(A;L)$ is a representation of the sub-adjacent Lie algebra $A^c$. By Lemma \ref{lem:den-ass},  $(A;L_{\succ},R_\prec)$ is a representation of the associative algebra $(A,\cdot)$. Moreover, \eqref{eq:pre-Poisson 1} implies that \eqref{eq:rep 1} holds, \eqref{eq:pre-Poisson 2} implies that \eqref{eq:rep 2} holds and \eqref{eq:pre-Poisson 3} implies that \eqref{eq:con-rep} holds.  Thus $(A;L_{\succ},R_\prec,L)$ is a representation of the noncommutative Poisson algebra $A^c$.\qed

\begin{cor}
  Let $(A,\succ,\prec,\ast)$ be a coherent noncommutative  pre-Poisson algebra. Then  the dual  $(A^*;-R_\prec^*, -L_\succ^*,L^*)$ is also a representation of  $A^c$, where $L_\succ^*,~R_\prec^*:A\longrightarrow \gl(V^*)$ are given by
$$
 \langle L_{\succ x}^*\alpha,y\rangle=-\langle \alpha,x\succ y\rangle,\quad\langle R_{\prec x}^*\alpha,y\rangle=-\langle \alpha,y\prec x\rangle,\quad \forall~x,y\in A,\alpha\in A^*.
$$
\end{cor}

\subsection{Rota-Baxter operators and $\huaO$-operators on noncommutative Poisson algebras }

A linear map $T:V\longrightarrow A$ is called an {\bf $\huaO$-operator} on an associative algebra $(A,\cdot)$ with respect to a representation  $(V;\huaL,\huaR)$ if $T$ satisfies
  \begin{equation}
    T(u)\cdot T(v)=T(\huaL_{T(u)}v+\huaR_{T(v)}u),\quad\forall~u,v\in V.
  \end{equation}
In particular, an $\huaO$-operator on an associative algebra $(A,\cdot)$ with respect to the regular representation is called a Rota-Baxter operator on $A$.
\begin{lem}{\rm(\cite{BGN2013A})}
  Let $(A,\cdot)$ be an associative algebra and $(V;\huaL,\huaR)$ a representation. Let $T:V\rightarrow A$ be an $\huaO$-operator on  $(A,\cdot)$ with respect to $(V;\huaL,\huaR)$. Then there exists a dendriform algebra structure on $V$ given by
$$
    u\succ v=\huaL_{T(u)}v,\quad u\prec v=\huaR_{T(v)}u, \quad\forall~u,v\in V.
$$
\end{lem}

A linear map $T:V\longrightarrow \g$ is called an {\bf $\huaO$-operator} (\cite{Kuper1})
on a Lie algebra $(\g,[-,-]_\g)$ with respect to a representation $(V;\rho)$ if $T$ satisfies
\begin{equation}
  [T(u), T(v)]_\g=T\Big(\rho(T(u))(v)-\rho(T(v))(u)\Big),\quad \forall~u,v\in V.
\end{equation}
In particular, an $\huaO$-operator on a Lie algebra $(\g,[-,-]_\g)$ with respect to the adjoint representation is called a Rota-Baxter operator on $\g$.

\begin{lem}{\rm(\cite{Bai2007})}
Let $T:V\to \g$ be an $\huaO$-operator  on a Lie algebra $(\g,[-,-]_\g)$ with respect to a representation $(V;\rho)$. Define a multiplication $\ast$ on $V$ by
\begin{equation}
  u\ast v=\rho(Tu)(v),\quad \forall~u,v\in V.
\end{equation}
Then $(V,\ast)$ is a pre-Lie algebra.
 \end{lem}



 Let  $(V;\huaL,\huaR,\rho)$ be a representation of a noncommutative Poisson algebra $(P,\cdot_P,\{-,-\}_P)$.
\begin{defi}
  \begin{itemize}\item[{\rm(i)}]A linear operator $T:V\longrightarrow
P$ is called  an {\bf $\huaO$-operator on  $P$} if $T$ is both an $\huaO$-operator on the associative algebra $(P,\cdot_P)$ and an $\huaO$-operator on the Lie algebra $(P,\{-,-\}_P)$;
\item[{\rm(ii)}] A linear operator $\huaB:P\longrightarrow
P$ is called a {\bf Rota-Baxter operator on $P$} if $\huaB$ is both a Rota-Baxter operator on the associative algebra $(P,\cdot_P)$ and a Rota-Baxter operator on the Lie algebra $(P,\{-,-\}_P)$.
\end{itemize}
\end{defi}

When $(P,\cdot_P,\{-,-\}_P)$ is a usual Poisson algebra, i.e.   $\cdot_P$ is commutative,  we recover the
notion of a  Rota-Baxter operator  on a  Poisson algebra
introduced by Aguiar in \cite{A2}.

It is obvious that $\huaB: P\longrightarrow P$ is a Rota-Baxter operator on $P$ if and only if $\huaB$ is an  $\huaO$-operator on $P$ with respect to the representation $(P;L,R,\ad)$.

\begin{ex}{\rm
    Let $\huaB$ be a Rota-Baxter operator on an associative algebra $(A,\cdot)$. Then $\huaB$ is a Rota-Baxter operator on the coherent noncommutative Poisson algebra $(A,\cdot,\{-,-\}_\hbar)$ given by Example \ref{ex:coherent poisson 2}.
    }
\end{ex}

\begin{ex}\label{ex:id}{\rm
   Let $(A,\succ,\prec,\ast)$ be a noncommutative pre-Poisson algebra.  Then the identity map $\id$ is an $\huaO$-operator on $A^c$ with respect  to the representation $(A;L_\succ,R_\prec,L)$.}
   \end{ex}

Obviously, we have

\begin{pro}
  Let $(P,\cdot_P,\Poisson{-,-}_P)$ be a coherent noncommutative Poisson algebra and $r\in \wedge^2 P$ a solution of the $\PYBE$. Then $r^\sharp:P^*\longrightarrow P$ is an $\huaO$-operator on $P$ with respect  to the representation $(P^*;-R^*,-L^*,\ad^*)$.
\end{pro}

 An $\huaO$-operator on a noncommutative Poisson algebra gives
a noncommutative  pre-Poisson algebra.

\begin{thm}
  Let $(P,\cdot_P,\{-,-\}_P)$ be a noncommutative Poisson algebra and $T:V\longrightarrow
P$ an $\huaO$-operator on  $P$ with respect to the representation $(V;\huaL,\huaR,\rho)$. Define new operations $\succ,~\prec$ and $\ast$ on $V$ by
$$ u\succ v=\huaL_{T(u)}v,\quad u\prec v=\huaR_{T(v)}u,\quad u\ast v=\rho(T(u))v.$$
Then $(V,\succ,\prec,\ast)$ is a noncommutative  pre-Poisson algebra and $T$ is a homomorphism from $V^c$ to $(P,\cdot_P,\{-,-\}_P)$.

 Furthermore, if the representation satisfies \eqref{eq:coherent condition}, then $(V,\succ,\prec,\ast)$ is a coherent noncommutative  pre-Poisson algebra.
\end{thm}
\pf First by the fact that $T$ is an $\huaO$-operator on the associative algebra $(P,\cdot_P)$ as well as an $\huaO$-operator on the Lie algebra $(P,\{-,-\}_P)$ with respect to the representations $(V;\huaL,\huaR)$ and $(V;\rho)$ respectively, we deduce that $(V,\succ,\prec)$ is a dendriform algebra and $(V,\ast)$ is a pre-Lie algebra.

Denote by $\{u,v\}_T:=u\ast v-v\ast u$, then by the fact $T([u,v]_T)=\{T(u),T(v)\}_P$ and \eqref{eq:rep 1},
\begin{eqnarray*}
  &&(u\ast v-v\ast u)\succ w-u\ast(v\succ w)+v\succ(u\ast w)\\
  &=&\{u,v\}_T\succ w-u\ast(v\succ w)+v\succ(u\ast w)\\
  &=&\huaL_{T\{u,v\}_T} w-\rho(T(u))\huaL_{T(v)}w+\huaL_{T(v)}\rho(T(u))w\\
  &=&\huaL_{\{T(u),T(v)\}_P} w-\rho(T(u))\huaL_{T(v)}w+\huaL_{T(v)}\rho(T(u))w=0,
\end{eqnarray*}
which implies that \eqref{eq:pre-Poisson 1} holds. Similarly, by \eqref{eq:rep 2}, we can show that \eqref{eq:pre-Poisson 2} also holds.

Denote by $u\cdot_T v:=u\succ v+v\prec u$, then by the fact $T(u\cdot_T v)=T(u)\cdot_P T(v)$ and \eqref{eq:con-rep},
\begin{eqnarray*}
  &&(u\succ v+u\prec v)\ast w-(u\ast w)\prec v-u\succ(v\ast w)\\
  &=&(u\cdot_T v)\ast w-(u\ast w)\prec v-u\succ(v\ast w)\\
  &=&\rho(T(u\cdot_T v))w-\huaR_{T(v)}\rho(T(u))w-\huaL_{T(u)}\rho(T(v))w\\
  &=&\rho(T(u)\cdot_P T(v))w-\huaR_{T(v)}\rho(T(u))w-\huaL_{T(u)}\rho(T(v))w=0,
\end{eqnarray*}
which implies that \eqref{eq:pre-Poisson 3} holds. Thus, $(V,\succ,\prec,\ast)$ is a noncommutative  pre-Poisson algebra. It is obvious that $T$ is a homomorphism from $V^c$ to $(P,\cdot_P,\{-,-\}_P)$.

If the representation satisfies \eqref{eq:coherent condition}, then we have
\begin{eqnarray*}
 &&(u\succ v+v\prec u)\ast w-u\ast(v\succ w)-v\ast(w\prec u)\\
  &=&\rho(T(\huaL_{T(u)}v+\huaR_{T(v)}u))w-\rho(T(u))(\huaL_{T(v)}w)-\rho(T(v))(\huaL_{T(u)}w)\\
  &=&\rho(T(u)\cdot_P T(v))w-\rho(T(u))(\huaL_{T(v)}w)-\rho(T(v))(\huaL_{T(u)}w)=0,
\end{eqnarray*}
which implies that $(V,\succ,\prec,\ast)$ is  coherent.\qed

\begin{cor}
Let $(P,\cdot_P,\{-,-\}_P)$ be a noncommutative Poisson algebra and $T:V\longrightarrow
P$ an $\huaO$-operator on   $P$ with respect to the representation $(V;\huaL,\huaR,\rho)$. Then $T(V)=\{T(v)\mid v\in V\}\subset P$ is a subalgebra of $P$ and there is an induced noncommutative pre-Poisson algebra structure on $T(V)$ given by
$$T(u)\succ T(v)=T(u\succ v),\quad T(u)\prec T(v)=T(u\prec v),\quad T(u)\ast T(v)=T(u\ast v)$$
for all $u,v\in V$.
\end{cor}

\begin{cor}
 Let $(P,\cdot_P,\{-,-\}_P)$ be a noncommutative Poisson algebra. There is a noncommutative pre-Poisson algebra structure on $P$ such that its sub-adjacent noncommutative Poisson algebra
is exactly  $(P,\cdot_P,\{-,-\}_P)$
  if and only if there exists an invertible $\huaO$-operator on $(P,\cdot_P,\{-,-\}_P)$.
\end{cor}

\pf If $T:V\longrightarrow
P$ is an invertible $\huaO$-operator on  $P$ with respect to the representation $(V;\huaL,\huaR,\rho)$, then the compatible noncommutative pre-Poisson algebra structure on $P$ is given by
$$x\succ y=T(\huaL_{x}T^{-1}(y)),\quad x\prec y=T(\huaR_{y}T^{-1}(x)),\quad x\ast y=T(\rho(x)(T^{-1}(y)))$$
for all $x,y\in P$.

Conversely, let $(P,\succ,\prec,\ast)$ be a noncommutative pre-Poisson algebra and $(P,\cdot_P,\{-,-\}_P)$ the sub-adjacent noncommutative Poisson algebra. Then the identity map $\id$ is an $\huaO$-operator on $P$ with respect to the representation $(P;L_{\succ},R_\prec,L)$. \qed

\begin{ex}{\rm
  Let $(P,\cdot_P,\{-,-\}_P)$ be a  noncommutative Poisson algebra and $\huaB:P\longrightarrow
P$ a Rota-Baxter operator. Define new operations on $P$ by
$$x\succ y=\huaB(x)\cdot_P y,\quad x\prec y=x\cdot_P \huaB(y),\quad x\ast y=\{\huaB(x),y\}_P.$$
Then $(P,\succ,\prec,\ast)$ is a noncommutative  pre-Poisson algebra and $\huaB$ is a homomorphism from the sub-adjacent noncommutative Poisson algebra $(P,\cdot_\huaB,\{-,-\}_\huaB)$ to $(P,\cdot_P,\{-,-\}_P)$, where $x\cdot_\huaB y=x\succ y+x\prec y$ and $\{x,y\}_\huaB=x\ast y-y\ast x$.
}
\end{ex}

\begin{ex}{\rm
  Let $(P,\cdot_P,\Poisson{-,-}_P)$ be a coherent noncommutative Poisson algebra and $r\in \wedge^2 P$ a solution of the $\PYBE$. Then $(P^*,\succ,\prec,\ast)$ is a coherent noncommutative  pre-Poisson algebra, where
  \begin{eqnarray*}
    \alpha\succ \beta=-R^*_{r^\sharp(\alpha)}\beta,\quad \alpha\prec \beta=-L^*_{r^\sharp(\beta)}\alpha,\quad \alpha\ast \beta=\ad^*_{r^\sharp(\alpha)}\beta,\quad\forall~\alpha,\beta\in P^*.
  \end{eqnarray*}
  }
\end{ex}

By Proposition \ref{pro:cocycle and symplectic}, we have
\begin{ex}{\rm
  Let $(P,\cdot_P,\Poisson{-,-}_P)$ be a coherent noncommutative Poisson algebra and $\omega\in \wedge^2 P^*$ non-degenerate. If $\omega$ is both a Connes cocycle on the associative algebra $(P,\cdot_P)$ and a symplectic structure on the Lie algebra $(P,\Poisson{-,-}_P)$, then $(P,\succ,\prec,\ast)$ is a noncommutative  pre-Poisson algebra, where $\succ,\prec$ and $\ast$ are determined by
  \begin{eqnarray*}
    \omega(x\succ y,z)=\omega(y,z\cdot_P x),\quad \omega(x\prec y,z)=\omega(x,y\cdot_P z),\quad \omega(x\ast y,z)=-\omega(y,\{x,z\}_P),\quad\forall~x,y,z\in P.
  \end{eqnarray*}
  }
\end{ex}

 \begin{thm}\label{thm:O-operator to Poisson bialgebra}
  Let $(V;\huaL,\huaR,\rho)$ be a representation of a coherent noncommutative Poisson algebra $(P,\cdot_P,\Poisson{-,-}_P)$  satisfying \eqref{eq:coherent condition}. Let $T:V\longrightarrow P$ be a linear map which is identified with an element in $(P\ltimes_{(-\huaR^*,-\huaL^*,\rho^*)} V^*)\otimes (P\ltimes_{(-\huaR^*,-\huaL^*,\rho^*)} V^* )$. Then $\overline{T}=T-\tau(T)$ is a skew-symmetric solution of the $\PYBE$ in $P\ltimes_{(-\huaR^*,-\huaL^*,\rho^*)} V^*$ if and only if $T$ is an $\huaO$-operator on $P$ with respect to the representation $(V;\huaL,\huaR,\rho)$, where $\tau$ is the exchange operator given by \eqref{eq:exchange operator}.
\end{thm}
\pf    By Proposition \ref{pro:semi-direct} and the fact that $(V;\huaL,\huaR,\rho)$ is a representation  satisfying \eqref{eq:coherent condition}, it follows that $P\ltimes_{(-\huaR^*,-\huaL^*,\rho^*)} V^*$ is a coherent noncommutative Poisson algebra.

By Corollary 3.10 in \cite{BGN2013A}, $\overline{T}=T-\tau(T)$ is a skew-symmetric solution of the associative Yang-Baxter equation in $P\ltimes_{(-\huaR^*,-\huaL^*)} V^*$ if and only if $T$ is an $\huaO$-operator on $(P,\cdot_P)$ with respect  to the representation $(\huaL,\huaR)$. By the conclusion in \cite{Bai2007},  $\overline{T}=T-\tau(T)$ is a skew-symmetric solution of the classical Yang-Baxter equation in $P\ltimes_{\rho^*} V^*$ if and only if $T$ is an $\huaO$-operator on $(P,\{\cdot,\cdot\}_P)$ with respect to the representation $\rho$. Thus, $\overline{T}=T-\tau(T)$ is a skew-symmetric solution of the Poisson Yang-Baxter equation if and only if $\overline{T}=T-\tau(T)$  is both a solution of the associative Yang-Baxter equation in $P\ltimes_{(-\huaR^*,-\huaL^*)} V^*$ and  the classical Yang-Baxter equation in $P\ltimes_{\rho^*} V^*$, which implies that $\overline{T}=T-\tau(T)$ is a skew-symmetric solution of the Poisson Yang-Baxter equation in $P\ltimes_{(-\huaR^*,-\huaL^*,\rho^*)} V^*$ if and only if $T$ is an $\huaO$-operator associated to the representation $(V;\huaL,\huaR,\rho)$.\qed\vspace{3mm}

   By Example \ref{ex:id} and Theorem \ref{thm:O-operator to Poisson bialgebra}, we get

   \begin{cor}
   Let $(A,\succ,\prec,\ast)$ be a coherent noncommutative pre-Poisson algebra. Then $r=\sum_i(e_i\otimes e^*_i-e^*_i\otimes e_i)$ is a skew-symmetric solution of the $\PYBE$ in the coherent noncommutative Poisson algebra $A^c\ltimes_{(-R_\prec^*, -L_\succ^*,L^*)} A^*$, where $\{e_1,\cdots,e_n\}$ is a basis of $A$ and $\{e^*_1,\cdots,e^*_n\}$ is the dual basis.
\end{cor}

\begin{ex}{\rm
 Let $(A,\succ,\prec)$ be a dendriform algebra. By Example \ref{ex:coherent pre-Poisson}, $(A,\succ,\prec,\ast_\hbar)$
is a coherent noncommutative pre-Poisson algebra, where
$x\ast_\hbar y=\hbar(x\succ y-y\prec x)$ for all $x,y\in A.$  Thus
$r=\sum_i(e_i\otimes e^*_i-e^*_i\otimes e_i)$ is a skew-symmetric
solution of the $\PYBE$ in the coherent noncommutative
Poisson algebra $A^c\ltimes_{(-R_\prec^*, -L_\succ^*,L^*)} A^*$,
where $\{e_1,\cdots,e_n\}$ is a basis of $A$ and
$\{e^*_1,\cdots,e^*_n\}$ is the dual basis.}
\end{ex}

\noindent
{\bf Acknowledgements. } This research is supported by NSFC (11922110, 11901501, 11931009).   C. Bai is also supported by the Fundamental Research Funds for the Central Universities and Nankai ZhiDe Foundation.

Jiefeng Liu\\
 School of Mathematics and Statistics, Northeast Normal University, Changchun 130024, China\\
 Email:liujf12@126.com

 Chengming Bai\\
 Chern Institute of Mathematics and LPMC, Nankai University,
Tianjin 300071, China \\ Email:  baicm@nankai.edu.cn

   Yunhe Sheng\\
Department of Mathematics, Jilin University,
 Changchun 130012,  China\\
 Email: shengyh@jlu.edu.cn

\end{document}